\documentclass[a4paper,12pt]{amsart}
\usepackage{amsmath,amsfonts,amssymb}
\usepackage{amsthm}
\usepackage{amstext,amscd,mathabx,mathdots}
\usepackage{graphics,graphicx,epsf}
\usepackage[usenames]{color}
\usepackage{hyperref}
\usepackage{tikz}
\usepackage{xcolor}
\usepackage{color}
\usepackage{float}
\usetikzlibrary{arrows}
\usepackage{lineno}
\usepackage{ragged2e}
\usepackage{calculator}

\theoremstyle{emsthmsl}

\newtheorem{theorem}{Theorem}[section]
\newtheorem{definition}[theorem]{Definition}
\newtheorem{proposition}[theorem]{Proposition}
\newtheorem{lemma}[theorem]{Lemma}

\newtheorem{remark}[theorem]{Remark}

\renewcommand{\proof}{{\noindent \bf Proof:\ }}

\numberwithin{equation}{section}

\title{Well-posedness for some third-order evolution differential equations: A semigroup approach} 

\author[F. D. M. Bezerra]{Flank D. M. Bezerra}
\address[F. D. M. Bezerra]{Departamento de Matem\'atica, Universidade Federal da Para\'iba, 58051-900 Jo\~ao Pessoa PB, Brazil.}
\email{flank.bezerra@academico.ufpb.br}

\author[A. N. Carvalho]{Alexandre N.
Carvalho$^\ddag$}\thanks{$^\ddag$Research partially supported by
CNPq \# 303929/2015-4 and by FAPESP \# 2003/10042-0, Brazil}
\address[A. N. Carvalho]{Departamento de Matem\'atica, Instituto de Ci\^encias
Mate\-m\'a\-ti\-cas e de Computa\c{c}\~{a}o, Universidade de S\~{a}o
Paulo-Campus de S\~{a}o Carlos, Caixa Postal 668, 13560-970 S\~{a}o
Carlos SP, Brazil.}
\email{andcarva@icmc.usp.br}

\author[L. A. Santos]{Lucas A. Santos}
\address[L. A. Santos]{Instituto Federal da Para\'iba, 58780-000 Itaporanga PB, Brazil.}
\email{lucas92mat@gmail.com}

\date{\today}

\begin{document}

\maketitle

\begin{abstract}

In this paper,  we discuss the well-posedness  of the Cauchy problem associated with the third-order evolution equation in time
$$
u_{ttt} +A u + \eta A^{\frac13} u_{tt} +\eta  A^{\frac23} u_t=f(u)
$$
where $\eta>0$, $X$ is a separable Hilbert space, $A:D(A)\subset X\to X$ is an unbounded sectorial operator with compact resolvent, and for some  $\lambda_0>0$ we have $\mbox{Re}\sigma(A)>\lambda_0$ and $f:D(A^{\frac13})\subset X\to X$ is a  nonlinear function  with suitable conditions of growth and regularity.


\end{abstract}

\tableofcontents

\section{Introduction} 

In this paper, we discuss the well-posedness of the Cauchy problem associated with the following  third-order evolution equation in time
\begin{equation}\label{Eq1}
u_{ttt} +A u + \eta A^{\frac13} u_{tt} + \eta A^{\frac23} u_t=f(u)
\end{equation}
where $\eta>0$, $X$ is a separable Hilbert space and $A:D(A)\subset X\to X$ is an unbounded sectorial operator with compact resolvent, and for some  $\lambda_0>0$ we have $\mbox{Re}\sigma(A)>\lambda_0$, that is,   $\mbox{Re}\lambda>\lambda_0$ for all $\lambda\in\sigma(A)$, where $\sigma(A)$ is the spectrum of $A$. This allows us to define the fractional power $A^{-\alpha}$ of order $\alpha\in(0,1)$ according to  \cite[Formula 4.6.9]{A} and  \cite[Theorem 1.4.2]{He}, as a closed linear operator on its domain  with inverse $A^{\alpha}$.

Denote by $X^{\alpha}=D(A^{\alpha})$ for $\alpha\in[0,1)$, taking $A^0:=I$ on $X^0:=X$ when $\alpha=0$. Recall that $X^{\alpha}$ is dense in $X$  for all $\alpha\in(0,1]$, for details see  \cite[Theorem 4.6.5]{A}. The fractional power space $X^\alpha$ endowed with the norm 
\[
\|\cdot\|_{X^\alpha}:=\|A^{\alpha} \cdot\|_X
\] 
is a Banach space. It is not difficult to show that $A^{\alpha}$ is the generator of a strongly continuous analytic
semigroup on $X$, that we will denote by $\{e^{-tA^{\alpha}}:  t\geqslant 0\}$, see  \cite{He} for any $\alpha\in[0,1]$. With this notation, we have $X^{-\alpha}=(X^\alpha)'$ for all $\alpha>0$, see  \cite{A}  for the characterization of the negative scale. 

Let $X_{-1}$ denote the extrapolation space of $X$ generated by $A$, and let $\{X_{-1}^\alpha:\alpha\geqslant0\}$ the fractional power scale generated by operator $A$ in $X_{-1}$, see \cite{A} and \cite{ACar} for more details.  

Here $f:D(A^{\frac13})\subset X\to X$ is  a  nonlinear function  with  suitable growth conditions and regularity in \eqref{Eq1} for different cases of $\eta>0$; namely, we consider:
\begin{itemize}
\item If $0<\eta< 1$, then we prove that the Cauchy problem defined by the linear equation associated with \eqref{Eq1} is ill-posed, consequently, the Cauchy problem defined by \eqref{Eq1} is ill-posed for any nonlinear function $f$, under the point of view of the theory of strongly continuous semigroups of bounded linear operators; 
\item If $\eta= 1$, then we assume that $f$ is twice continuously Fr\'echet differentiable and  Lipschitz continuous on bounded sets;
\item If $\eta>1$,  then we assume that $f$ is an $\epsilon$-regular map relative to the pair $(X^{\frac13},X)$ for $\epsilon\geqslant0$; that is, there exist  constants $c>0$, $\rho>1$, $\gamma(\epsilon)$ with $\rho\epsilon\leqslant \gamma(\epsilon)<\frac{1}{3}$ such that $f:X^{\frac{1}{3}+\epsilon}\to X^{\gamma(\epsilon)}$ and  
\begin{equation}\label{asaE4}
\|f(\phi_1)-f(\phi_2)\|_{X^{\gamma(\epsilon)}}\leqslant c\|\phi_1-\phi_2\|_{X^{\frac{1}{3}+\epsilon}}(1+\|\phi_1\|_{X^{\frac{1}{3}+\epsilon}}^{\rho-1}+\|\phi_2\|_{X^{\frac{1}{3}+\epsilon}}^{\rho-1}),
\end{equation}
for any $\phi_1,\phi_2\in X^{\frac{1}{3}+\epsilon}$, see \cite[Definition 2]{ACar}, \cite{CCCC} and \cite{CC_2} for more details.
\end{itemize}

For a better understanding of the $\epsilon$-regular map relative to the pair $(X^{\frac13},X)$ for $\epsilon\geqslant0$, we construct the following diagram.
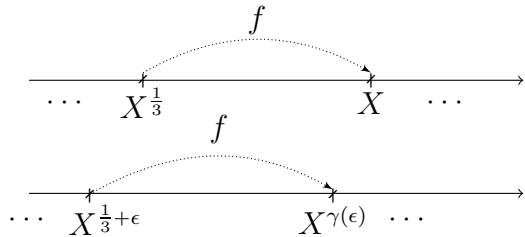
\begin{figure}[h]\label{Figureee}
\begin{tikzpicture}
\draw[->] (-3.5,0) -- (3,0);
\node at (1,0.01) {$\nshortmid$};
\node at (1,-0.3) {$X$};
\node at (-2,0.01) {$\nshortmid$};
\node at (-2,-0.3) {$X^\frac13$};
\node at (-0.5,0.8) {$f$};
\draw[densely dotted,-latex'] (-2,0.1) to [out=30,in=150] (1,0.1);
\draw[->] (-3.5,-1.5) -- (3,-1.5);
\node at (0.5,-1.51) {$\nshortmid$};
\node at (0.5,-1.85) {$X^{\gamma(\epsilon)}$};
\node at (-2.7,-1.51) {$\nshortmid$};
\node at (-2.5,-1.85) {$X^{\frac{1}{3}+\epsilon}$};
\node at (-1,-0.65) {$f$};
\draw[densely dotted,-latex'] (-2.7,-1.51) to [out=30,in=150] (0.5,-1.45);
\node at (-3,-0.3) {$\cdots$};
\node at (-3.5,-1.85) {$\cdots$};
\node at (1.5,-1.85) {$\cdots$};
\node at (2,-0.3) {$\cdots$};
\end{tikzpicture}
\caption{$X^{\frac{1}{3}+\epsilon}\subset X^{\frac13}$ and $X^{\gamma(\epsilon)}\subset X$}
\end{figure}

The evolution equations of third order in time have been studied extensively in the Hilbert setting and much progress has been achieved. In \cite{ALM,CLC,CLR2,KLM,KLP,LW,MMT,PS-H} and \cite{PS-M}, and references therein, the Moore-Gibson-Thompson (MGT) equation is studied in different contexts and results of existence, stability and regularity of solutions are obtained by the spectral theory of the MGT operator. In \cite{HOF} the abstract linear equations of third order in time is analyzed and results on (non)existence of solution are obtained. In \cite{BS} the abstract linear equations of third order in time is analyzed and results on (non) existence, stability and regularity of solution are obtained via theory of fractional powers of closed and densely defined operators.

The article is organized in the following way. In Section \ref{Sec: SpecProp} we present general facts on spectral behavior of the our problem.   In Section \ref{S0} we consider the  case $0\leqslant\eta<1$  and we obtain the result that shows that the problem \eqref{Eq1} is ill-posed under the point of view of the theory of strongly continuous semigroups of bounded linear operators. In Section \ref{S3}  we consider the  case $\eta>1$ and we obtain a result of existence, stability and regularity of solutions for \eqref{Eq1}.   In Section \ref{S1}   we consider the  case $\eta=1$ via theory of strongly continuous groups of bounded linear operators.  Finally, in Section \ref{FinalSec}   we explore our results on the Moore-Gibson-Thompson-type equations, according to the references \cite{ALM,CLC,CLR2,KLM,KLP,LW,MMT,PS-H} and \cite{PS-M}.

\section{Functional Framework}\label{Sec: SpecProp}

We first introduce some notations, we consider $Z=X^{\frac23}\times X^{\frac13}\times X$ endowed with the norm given by
\[
\Big\|\Big[\begin{smallmatrix} u\\ v\\ w \end{smallmatrix}\Big]\Big\|^2_{Z}=\|u\|^2_{X^{\frac23}}+\|v\|^2_{X^{\frac13}}+\|w\|^2_X,\ \forall \Big[\begin{smallmatrix} u\\ v\\ w \end{smallmatrix}\Big]\in Z.
\] 

\subsection{Natural setting}

We can rewrite the initial value problem associated with equation \eqref{Eq1} as the Cauchy problem in $Z$
\begin{equation}\label{Eq4aaa}
\dfrac{d}{dt}\Big[\begin{smallmatrix} u\\ v\\ w \end{smallmatrix}\Big]+ \mathbb{A}_{(\eta)}\Big[\begin{smallmatrix} u\\ v\\ w \end{smallmatrix}\Big]=F(\Big[\begin{smallmatrix} u\\ v\\ w \end{smallmatrix}\Big]),\ t>0,
\end{equation}
and
\begin{equation}\label{Eq4bbb}
\Big[\begin{smallmatrix} u\\ v\\ w \end{smallmatrix}\Big](0)=\Big[\begin{smallmatrix} u_0\\ v_0\\ w_0 \end{smallmatrix}\Big],
\end{equation}
where $v=u_t$ and $w=v_t$ and the unbounded linear operator $\mathbb{A}_{(\eta)}:D(\mathbb{A}_{(\eta)})\subset Z\to Z$ is defined by
\begin{equation}\label{definicaoOp01}
D(\mathbb{A}_{(\eta)})=X^1\times X^{\frac23}\times X^{\frac13}
\end{equation}
and
\begin{equation}\label{definicaoOp02}
\mathbb{A}_{(\eta)}\Big[\begin{smallmatrix} u\\ v\\ w \end{smallmatrix}\Big]:= \Big[\begin{smallmatrix} 0 & -I & 0 \\ 0 & 0 & -I\\ A & \eta A^{\frac23} & \eta A^{\frac13} \end{smallmatrix}\Big] \Big[\begin{smallmatrix} u\\ v\\ w \end{smallmatrix}\Big]=\Big[\begin{smallmatrix} -v\\ -w\\ Au+\eta A^{\frac23}v+\eta A^{\frac13}w \end{smallmatrix}\Big],\ \forall \Big[\begin{smallmatrix} u\\ v\\ w \end{smallmatrix}\Big]\in X^1\times X^{\frac23}\times X^{\frac13}.
\end{equation}

The nonlinearity $F$ given by 
\begin{equation}\label{DefF}
F(\Big[\begin{smallmatrix} u\\ v\\ w \end{smallmatrix}\Big])=\Big[\begin{smallmatrix} 0\\ 0\\ f(u) \end{smallmatrix}\Big],
\end{equation}
where $f:D(A^{\frac13})\subset X\to X$ is a Lipschitz continuous function on bounded sets.

From now on, we denote 
\[
Z^1=D(\mathbb{A}_{(\eta)})=X^1\times X^{\frac23}\times X^{\frac13}.
\]

We also consider the following notion of mild solution for \eqref{Eq4aaa}-\eqref{Eq4bbb}. Given $\left[\begin{smallmatrix} u_0\\ v_0\\ w_0\end{smallmatrix}\right]\in Z$ we say that $\left[\begin{smallmatrix} u\\ v\\ w\end{smallmatrix}\right]$ is a local  mild solution of \eqref{Eq4aaa}-\eqref{Eq4bbb} provided $\left[\begin{smallmatrix} u\\ v\\ w\end{smallmatrix}\right]\in C([0,\tau_{u_{0}, v_{0}, w_0}), Z)$, $f(u)\in C([0,\tau_{u_{0}, v_{0},w_0}), X)$ and, for $t\in(0,\tau_{u_{0}, v_{0},w_0})$, $\left[\begin{smallmatrix} u\\ v\\ w\end{smallmatrix}\right]$ satisfies the integral equation
\begin{equation}\label{Newrwgery}
\left[\begin{smallmatrix} u(t)\\ v(t)\\ w(t)\end{smallmatrix}\right]=e^{-\mathbb{A}_{(\eta)} t}\left[\begin{smallmatrix} u_0\\ v_0\\ w_0\end{smallmatrix}\right]
+\int_0^t e^{-\mathbb{A}_{(\eta)}(t-s)}\left[\begin{smallmatrix} 0\\ 0\\ f(u(s))\end{smallmatrix}\right]ds.
\end{equation}
for some $\tau_{u_{0}, v_{0},w_0}>0$.

\begin{lemma}\label{Mle}
Let $\mathbb{A}_{(\eta)}$ be the  unbounded linear operator  defined in \eqref{definicaoOp01}-\eqref{definicaoOp02}. The following conditions hold.
\begin{itemize}
\item[$i)$] The linear operator $\mathbb{A}_{(\eta)}$ is closed and densely defined;
\item[$ii)$] Zero belongs to the resolvent set $\rho(\mathbb{A}_{(\eta)})$; namely, the resolvent operator of $\mathbb{A}_{(\eta)}$ is the bounded linear operator $\mathbb{A}_{(\eta)}^{-1}:Z\to Z$ given by
\[
\mathbb{A}_{(\eta)}^{-1}=\Big[\begin{smallmatrix} \eta A^{-\frac13} & \eta A^{-\frac23} & A^{-1}\\ -I& 0& 0 \\ 0 & -I & 0 \end{smallmatrix}\Big].
\]
Moreover, $\mathbb{A}_{(\eta)}$ has compact resolvent;
\item[$iii)$] The spectrum set of $-\mathbb{A}_{(\eta)}$, $\sigma(-\mathbb{A}_{(\eta)})$, is given by
\[
\sigma(-\mathbb{A}_{(\eta)})=\{\lambda \in \mathbb{C}: \lambda \in \sigma(-A^{\frac13})\}\cup\{z_\eta\lambda \in \mathbb{C}: \lambda \in \sigma(-A^{\frac13})\}\cup\{\overline{z_\eta}\lambda \in \mathbb{C}: \lambda \in \sigma(-A^{\frac13})\},
\]
where $\sigma(-A^{\frac13})$ denote the spectrum set of  $-A^{\frac13}$ and
\[
z_\eta=\dfrac{1}{2}\Big[(\eta-1)-i\sqrt{3+2\eta-\eta^2}\Big]
\]
and
\[
\overline{z_\eta}=\dfrac{1}{2}\Big[(\eta-1)+i\sqrt{3+2\eta-\eta^2}\Big].
\]
\end{itemize}
\end{lemma}

\proof
To prove part  $i)$ we take a sequence $\Big(\Big[\begin{smallmatrix} u_n\\v_n\\ w_n \end{smallmatrix}\Big],\mathbb{A}_{(\eta)}\Big[\begin{smallmatrix} u_n\\v_n\\ w_n \end{smallmatrix}\Big] \Big)$ is the graph of $\mathbb{A}_{(\eta)}$, which converges in $ Z$ to $\Big(\Big[\begin{smallmatrix} \phi\\ \varphi\\ \psi \end{smallmatrix}\Big],\Big[\begin{smallmatrix} \nu\\ \mu\\ \zeta \end{smallmatrix}\Big]   \Big)$. From this we easily conclude that $\nu=-\varphi$, $\mu=-\psi$, and 
\[
u_n\to \phi\ \mbox{in}\ X^{\frac23}\Longrightarrow A^{\frac23}u_n\to A^{\frac23}\phi\ \mbox{in}\ X,
\] 
\[
v_n\to \varphi\ \mbox{in}\ X^{\frac13}\Longrightarrow A^{\frac13}v_n\to A^{\frac13}\varphi\ \mbox{in}\ X,
\]
and
\[
Au_n+\eta A^{\frac23}v_n+\eta A^{\frac13}w_n\to \zeta\ \mbox{in}\ X.
\]

Since 
\[
 A^{\frac23}u_n+\eta A^{\frac13}v_n+\eta w_n\to A^{\frac23}\phi+\eta A^{\frac13}\varphi+\eta \psi\ \mbox{in}\ X,
\]
the  closedness of $A^{\frac13}$ imples that $A^{\frac23}\phi+\eta A^{\frac13}\varphi+\eta \psi\in X^{\frac13}$ and $A\phi+\eta A^{\frac23}\varphi+\eta A^{\frac13}\psi=\zeta$, and we obtain the result; that is, $\Big[\begin{smallmatrix} \phi\\ \varphi\\ \psi \end{smallmatrix}\Big]\in D(\mathbb{A}_{(\eta)})$ and $\mathbb{A}_{(\eta)}\Big[\begin{smallmatrix} \phi\\ \varphi\\ \psi \end{smallmatrix}\Big]=\Big[\begin{smallmatrix} \nu\\ \mu\\ \zeta \end{smallmatrix}\Big]$.

For the proof of $ii)$ the result easily follows from $\mathbb{A}_{(\eta)}^{-1}:Z\to Z$ to be given by
\[
\mathbb{A}_{(\eta)}^{-1}=\Big[\begin{smallmatrix} \eta A^{-\frac13} & \eta A^{-\frac23} & A^{-1}\\ -I& 0& 0 \\ 0 & -I & 0 \end{smallmatrix}\Big]
\]
which takes bounded subsets of $X^{\frac23}\times X^{\frac13}\times X$ into bounded subsets of $X^{1}\times X^{\frac23}\times X^{\frac13}$. The latter space is compactly embedded in $Z$ because the inclusions $X^{\beta}\subset X^{\gamma}$ are compact, for $\beta>\gamma\geq 0$, provided $A$ has compact resolvent.   
Finally, after considering the eigenvalue problem for the operator  $-\mathbb{A}_{(\eta)}$, 
\[
-\mathbb{A}_{(\eta)}{\bf u} = \lambda{\bf u},
\]
and after straightforward calculations, part $iii)$ follows immediately from the fact that $\mathbb{A}_{(\eta)}$ has compact resolvent.
\qed

\subsection{New approach}\label{Subsection2.2}

Let \eqref{Eq1} be the original equation of third-order.  Note that using the change variable 
\[
v=u_t + A^{\frac13}u,
\]
and the function
\[
w=v_t,
\]
we can rewrite the  \eqref{Eq1} as follows, a first-order evolution equation in time for $w$
\begin{equation}\label{E00q1Linear3}
w_t  + (\eta-1) A^{\frac13} w+  A^{\frac23} v =f(u).
\end{equation}

The initial value problem associated with equation \eqref{E00q1Linear3} as the Cauchy problem in $Z$
\begin{equation}\label{Eq4aadfdfdfabb}
\dfrac{d}{dt}\Big[\begin{smallmatrix} u\\ v\\ w \end{smallmatrix}\Big]+ \mathbb{B}_{(\eta)}\Big[\begin{smallmatrix} u\\ v\\ w \end{smallmatrix}\Big]=F(\Big[\begin{smallmatrix} u\\ v\\ w \end{smallmatrix}\Big]),\ t>0,
\end{equation}
and
\begin{equation}\label{Eq4bbbfdfdbbbb}
\Big[\begin{smallmatrix} u\\ v\\ w \end{smallmatrix}\Big](0)=\Big[\begin{smallmatrix} u_0\\ v_0\\ w_0 \end{smallmatrix}\Big],
\end{equation}
where $v=u_t$ and $w=v_t$ and the unbounded linear operator $\mathbb{B}_{(\eta)}:D(\mathbb{B}_{(\eta)})\subset Z\to Z$ is defined by
\begin{equation}\label{Eq4aaabb}
D(\mathbb{B}_{(\eta)})=X^1\times X^{\frac23}\times X^{\frac13},
\end{equation}
and
\begin{equation}\label{Eq4bbbbbbb}
\mathbb{B}_{(\eta)}\Big[\begin{smallmatrix} u\\ v\\ w \end{smallmatrix}\Big]:= \Big[\begin{smallmatrix} A^{\frac13} & -I & 0 \\ 0 & 0 & - I\\ 0 &  A^{\frac23} & (\eta-1) A^{\frac13} \end{smallmatrix}\Big] \Big[\begin{smallmatrix} u\\ v\\ w \end{smallmatrix}\Big]=\Big[\begin{smallmatrix} A^{\frac13}u-v\\ -w\\  A^{\frac23}v+(\eta-1) A^{\frac13}w \end{smallmatrix}\Big],\ \forall \Big[\begin{smallmatrix} u\\ v\\ w \end{smallmatrix}\Big]\in X^1\times X^{\frac23}\times X^{\frac13}.
\end{equation}

The nonlinearity $F$ given by \eqref{DefF}.

\begin{lemma}\label{MleNewddsd}
The following conditions hold.
\begin{itemize}
\item[$i)$] The linear operator $\mathbb{B}_{(\eta)}$ is closed and densely defined;
\item[$ii)$] Zero belongs to the resolvent set $\rho(\mathbb{B}_{(\eta)})$; namely, the resolvent operator of $\mathbb{B}_{(\eta)}$ is the bounded linear operator $\mathbb{B}_{(\eta)}^{-1}:Z\to Z$ given by
\[
\mathbb{B}_{(\eta)}^{-1}=\Big[\begin{smallmatrix}  A^{-\frac13} & (\eta-1) A^{-\frac23} & A^{-1}\\ 0& (\eta-1) A^{-\frac13}& A^{-\frac23}  \\ 0 & -I & 0 \end{smallmatrix}\Big].
\]
Moreover, $\mathbb{B}_{(\eta)}$ has compact resolvent;
\item[$iii)$] The spectrum set of $-\mathbb{B}_{(\eta)}$, $\sigma(-\mathbb{B}_{(\eta)})$, is given by
\[
\sigma(-\mathbb{B}_{(\eta)})=\{\lambda \in \mathbb{C}: \lambda \in \sigma(-A^{\frac13})\}\cup\{c_\eta\lambda \in \mathbb{C}: \lambda \in \sigma(-A^{\frac13})\}\cup\{d_\eta\lambda \in \mathbb{C}: \lambda \in \sigma(-A^{\frac13})\},
\]
where $\sigma(-A^{\frac13})$ denote the spectrum set of  $-A^{\frac13}$ and
\[
c_\eta=\dfrac{1}{2}\Big[(\eta-1)+\sqrt{\eta^2-2\eta-3}\Big]
\]
and
\[
d_\eta=\dfrac{1}{2}\Big[(\eta-1)-\sqrt{\eta^2-2\eta-3}\Big],
\]
\end{itemize}
\end{lemma}

\proof
To prove part  $i)$ we take a sequence $\Big(\Big[\begin{smallmatrix} u_n\\v_n\\ w_n \end{smallmatrix}\Big],\mathbb{B}_{(\eta)}\Big[\begin{smallmatrix} u_n\\v_n\\ w_n \end{smallmatrix}\Big] \Big)$ is the graph of $\mathbb{B}_{(\eta)}$, which converges in $ Z$ to $\Big(\Big[\begin{smallmatrix} \phi\\ \varphi\\ \psi \end{smallmatrix}\Big],\Big[\begin{smallmatrix} \nu\\ \mu\\ \zeta \end{smallmatrix}\Big]   \Big)$. Then we have 

\[
u_n\to \phi\ \mbox{in}\ X^{\frac23}
\] 
\[
v_n\to \varphi\ \mbox{in}\ X^{\frac13}
\]
\[
w_n\to \psi\ \mbox{in}\ X 
\]

and

\[
A^\frac13 u_n -v_n \to \nu\ \mbox{in}\ X^{\frac23}
\] 
\[
-w_n\to \mu\ \mbox{in}\ X^{\frac13}
\]
\[
A^\frac23 v_n+(\eta-1)A^\frac13 w_n\to \zeta\ \mbox{in}\ X \]

From this, we easily conclude that

$$\psi=-\mu\in X^\frac13.$$

From the closedness of $A^\frac13$, we obtain

$$A^\frac13\varphi+(\eta-1)\psi \in X^\frac13\ \  
\text{and}\ \  A^\frac13(A^\frac13\varphi+(\eta-1)\psi)=\zeta,$$
that is, 
$$\varphi\in X^\frac23  \  \text{and}\ \ A^\frac23\varphi+(\eta-1)A^\frac13\psi=\zeta$$

From the closedness of $A^\frac13$, we also obtain
$$A^\frac23 \phi-A^\frac13 \varphi \in X^\frac13 \ \ \text{and}\ \ A^\frac13(A^\frac23 \phi-A^\frac13 \varphi)=A^\frac23\nu $$
that is, 
$$\phi\in X^1 \ \ \text{and}\ \ A^\frac13\phi-\varphi=\nu $$

From this, we conclude that $\Big[\begin{smallmatrix} \phi\\ \varphi\\ \psi \end{smallmatrix}\Big]\in D(\mathbb{B}_{(\eta)})$ and $\mathbb{B}_{(\eta)}\Big[\begin{smallmatrix} \phi\\ \varphi\\ \psi \end{smallmatrix}\Big]=\Big[\begin{smallmatrix} \nu\\ \mu\\ \zeta \end{smallmatrix}\Big]$.

For the proof of $ii)$ the result easily follows from $\mathbb{B}_{(\eta)}^{-1}:Z\to Z$ to be given by
\[
\mathbb{B}_{(\eta)}^{-1}=\Big[\begin{smallmatrix} \eta A^{-\frac13} & (1-\eta) A^{-\frac23} & A^{-1}\\ 0& (1-\eta) A^{-\frac13}& A^{-\frac23}  \\ 0 & -I & 0 \end{smallmatrix}\Big].
\]

Finally, after considering the eigenvalue problem for the operator  $-\mathbb{B}_{(\eta)}$, 
\[
-\mathbb{B}_{(\eta)}{\bf u} = \lambda{\bf u},
\]
and after straightforward calculations, part $iii)$ follows immediately from the fact that $\mathbb{B}_{(\eta)}$ has compact resolvent.
\qed

\section{Ill-posed problems}\label{S0}

In this section we consider the case $0\leq\eta<1$. We prove that the Cauchy problem defined by the linear equation associated with \eqref{Eq1} is ill-posed in $Z$, consequently, the Cauchy problem defined by \eqref{Eq1} is ill-posed for any nonlinear function $f$ in $Z$, under the point of view of the theory of strongly continuous semigroups of bounded linear operators.

\begin{lemma}\label{lemma1}
Let $\mathbb{A}_{(\eta)}$ be the  unbounded linear operator  defined in \eqref{definicaoOp01}-\eqref{definicaoOp02}.  If $0\leq \eta < 1$, then the unbounded linear operator $-\mathbb{A}_{(\eta)}$ with $\mathbb{A}_{(\eta)}:D(\mathbb{A}_{(\eta)})\subset Z\to Z$  is not the infinitesimal generator of a strongly continuous semigroup in $Z$.
\end{lemma} 

\proof
If $-\mathbb{A}_{(\eta)}$ generates a strongly continuous semigroup $\{e^{-\mathbb{A}_{(\eta)} t}:t\geqslant0\}$ in $Z$, it follows from \cite{A} that there exist  constants $\omega({\eta)}\geq 0$ and $M(\eta)\geq 1$ such that

\begin{equation}
\|e^{-\mathbb{A}_{(\eta)} t}\|_{\mathcal{L}(Y)}\leq M(\eta)e^{\omega(\eta) t}\ \ \text{\ for \ } t\geqslant0.
\end{equation}
Moreover, from \cite{A} we have

\begin{equation}\label{semiplan}
\{\lambda \in \mathbb{C}: Re\lambda > \omega(\eta)\} \subset \rho(-\mathbb{A}_{(\eta)}) 
\end{equation}
where $\rho(-\mathbb{A}_{(\eta)})$ denotes the resolvent set of the operator $-\mathbb{A}_{(\eta)}$. From Lemma \ref{Mle}(iii) we have
\[
\sigma(-\mathbb{A}_{(\eta)})=\{\lambda \in \mathbb{C}: \lambda \in \sigma(-A^{\frac13})\}\cup\{z_\eta\lambda \in \mathbb{C}: \lambda \in \sigma(-A^{\frac13})\}\cup\{\overline{z_\eta}\lambda \in \mathbb{C}: \lambda \in \sigma(-A^{\frac13})\},
\]
where
\[
z_\eta=\dfrac{1}{2}\Big((\eta-1)+i\sqrt{3+2\eta-\eta^2}\Big)
\]
and
\[
\overline{z_\eta}=\dfrac{1}{2}\Big((\eta-1)-i\sqrt{3+2\eta-\eta^2}\Big),
\]
where $3+2\eta-\eta^2>0$ for any  $0\leq \eta < 1$, see Figure \ref{figadadaf1}.

	\begin{figure}[H]
		\begin{center}
			\begin{tikzpicture}
			\draw[-stealth'] (-5,0) -- (4,0) node[below] {$\scriptstyle {\rm Re}$};
			\draw[-stealth'] (0,-3.15) -- (0,3.15) node[above] {\color{black}$\scriptstyle{\rm Im}$};
			\draw[->] (4mm,0mm) arc (0:40:5mm);
			\node at (0.9,0.4) {$\theta_\eta$};
			\draw (0,0) -- (-3, 0);
\foreach \i in {0.1,0.14,...,1.2}{\EXP[\i]{3}{\sol}\fill [color=blue] (-3*\sol,0) circle (1.5pt);}
\node[blue] at (-5,0.4) {$\sigma(-A^{\frac13})$};
\foreach \i in {0.1,0.14,...,1.2}{\EXP[\i]{3}{\sol}\fill [color=blue] (1.3*\sol,2*\sol) circle (1.5pt);}
\node[blue] at (3.2,3) {$z_\eta\sigma(-A^{\frac13})$};
\foreach \i in {0.1,0.14,...,1.2}{\EXP[\i]{3}{\sol}\fill [color=blue] (1.3*\sol,-2*\sol) circle (1.5pt);}
\node[blue] at (3.2,-3) {$\overline{z_\eta}\sigma(-A^{\frac13})$};
			\end{tikzpicture}
		\end{center}
\caption{Semi-lines contained the eigenvalues of $-\mathbb{A}_{(\eta)}$ and $0<\theta_\eta<\dfrac{\pi}{2}$, $0\leqslant \eta < 1$.}\label{figadadaf1}
	\end{figure}
	
Since $\sigma(-A^{\frac13})=\{-\mu_n: n\in \mathbb{N}\}$ with $\mu_n\in \sigma(A^{\frac13})$ for each $n\in \mathbb{N}$ and $\mu_n\to \infty$ as $n\to \infty$ and $0\leq\eta<1$, we conclude that
\[
\sigma(-\mathbb{A}_{(\eta)})\cap \{\lambda \in \mathbb{C}:Re \lambda>\omega(\eta)\}\neq \emptyset
\]
This contradicts \eqref{semiplan} and therefore $-\mathbb{A}_{(\eta)}$ can not be the infinitesimal generator of a strongly continuous semigroup in $Z$.
\qed

\begin{theorem}\label{MThe}
If $0\leqslant \eta < 1$, then the Cauchy problem
\[
\dfrac{d}{dt}\Big[\begin{smallmatrix} u\\ v\\ w \end{smallmatrix}\Big]+ \mathbb{A}_{(\eta)}\Big[\begin{smallmatrix} u\\ v\\ w \end{smallmatrix}\Big]=\Big[\begin{smallmatrix} 0\\ 0\\ 0 \end{smallmatrix}\Big],\ t>0,
\]
and
\[
\Big[\begin{smallmatrix} u\\ v\\ w \end{smallmatrix}\Big](0)=\Big[\begin{smallmatrix} u_0\\ v_0\\ w_0 \end{smallmatrix}\Big],
\]
where $v=u_t$ and $w=v_t$ is ill-posed on $Z$.
\end{theorem}

\proof
The result easily follows from Lemma \ref{lemma1}.
\qed

\begin{theorem}
The Cauchy problem \eqref{Eq4aaa}-\eqref{Eq4bbb}  is ill-posed on $Z$. More precisely, let $\left[\begin{smallmatrix} u_0\\ v_0\\ w_0\end{smallmatrix}\right]\in Z$  does not exist $\left[\begin{smallmatrix} u\\ v\\ w\end{smallmatrix}\right]\in C([0,\tau_{u_{0}, v_{0}, w_0}), Z)$ with  $f(u)\in C([0,\tau_{u_{0}, v_{0},w_0}), X)$ such that \eqref{Newrwgery} holds, for any $\tau_{u_{0}, v_{0},w_0}>0$.
\end{theorem}

\proof
The result easily follows from Theorem \ref{MThe}.
\qed

\begin{remark}
If $\eta=0$ then thanks to the results in \cite{BS} the Cauchy linear problem associated with \eqref{Eq4aaa} is ill-posed in $Z$, and therefore, the Cauchy  problem  \eqref{Eq4aaa} is ill-posed in $Z$.
\end{remark}

On the new approach present in Subsection \ref{Subsection2.2} we have

\begin{theorem}
Let $0\leqslant\eta<1$ and let $\mathbb{B}_{(\eta)}$ be the unbounded linear operator defined in \eqref{Eq4aaabb}-\eqref{Eq4bbbbbbb}. Then the problem \eqref{Eq4aadfdfdfabb}-\eqref{Eq4bbbfdfdbbbb} is ill-posed in the sense that it does not generate a strongly continuous semigroup of bounded linear operators on the state space $Z$.
\end{theorem}

\proof
If $-\mathbb{B}_{(\eta)}$ generates a strongly continuous semigroup $\{e^{-\mathbb{B}_{(\eta)} t}:t\geqslant0\}$ on $Z$, it follows from  Pazy \cite[Theorem 1.2.2]{P} that there exist  constants $\omega\geq 0$ and $M\geq 1$ such that
\begin{equation}
\|e^{-\mathbb{B}_{(\eta)} t}\|_{\mathcal{L}(Y)}\leq Me^{\omega t}\ \ \text{\ for \ } 0\leq t< \infty.
\end{equation}
Moreover, from Pazy \cite[Remark 1.5.4]{P} we have
\begin{equation}\label{semiplan}
\{\lambda \in \mathbb{C}: Re\lambda > \omega\} \subset \rho(-\mathbb{B}_{(\eta)}). 
\end{equation}
where $\rho(-\mathbb{B}_{(\eta)})$ denotes the resolvent set of the operator $-\mathbb{B}_{(\eta)}$.

From Lemma \ref{MleNewddsd} we can consider a sequence $(\lambda_k z_{\eta})_{k}\in\sigma(-\mathbb{B}_{(\eta)})$, for $k=1,2,3,\dots$, with $\lambda_k\in\sigma(-A^{\frac13})$ and $|\lambda_k|\to \infty$ as $k\to \infty$. Note that 

$$arg(\lambda_k z_{\eta})=\arctan\left(\frac{2\sqrt{-\eta^2+2\eta+3}}{1-\eta}\right)$$
and since $0\leq\eta<1$, we have

$$0<arg(\lambda_k z_{\eta})<\pi/2$$
for every $k\geq 1$ and $|\lambda_k z_{\eta}| \to\infty$ as $k\to \infty$. Then we conclude that
\[
\sigma(-\mathbb{B}_{(\eta)})\cap \{\lambda \in \mathbb{C}:Re \lambda>\omega\}\neq \emptyset.
\]

This contradicts \eqref{semiplan} and therefore $-\mathbb{B}_{(\eta)}$ can not be the infinitesimal generator of a strongly continuous semigroup on $Z$.
\qed

\section{Parabolic differential equations}\label{S3}

In this section we consider the case $\eta>1$. Namely, thanks to the Lemma \ref{MleNewddsd} we have the following illustration of the eigenvalues of $-\mathbb{B}_{(\eta)}$.

	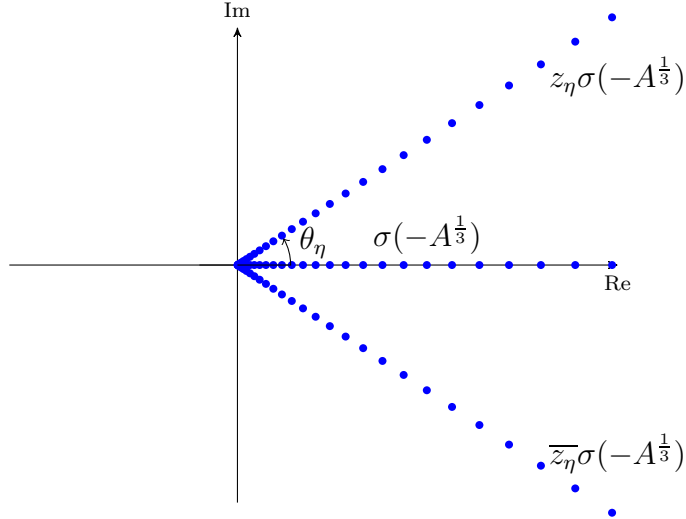
\begin{figure}[H]
		\begin{center}
			\begin{tikzpicture}
			\draw[-stealth'] (-0.5,0) -- (5,0) node[below] {$\scriptstyle {\rm Re}$};
			\draw[-stealth'] (0,-3.15) -- (0,3.15) node[above] {\color{black}$\scriptstyle{\rm Im}$};
			\draw (0,0) -- (-3, 0);
\foreach \i in {0.1,0.14,...,1.2}{\EXP[\i]{3}{\sol}\fill [color=blue] (3*\sol,0) circle (1.5pt);}
\foreach \i in {0.1,0.14,...,1.2}{\EXP[\i]{3}{\sol}\fill [color=blue] (3*\sol,2*\sol) circle (1.5pt);}
\foreach \i in {0.1,0.14,...,1.2}{\EXP[\i]{3}{\sol}\fill [color=blue] (3*\sol,-2*\sol) circle (1.5pt);}
\draw[->] (0,0) -- (7mm,0mm) arc (0:30:7mm);
\node at (1,0.3) {$\theta_\eta$};
\node at (2.5,0.4) {$\sigma(-A^{\frac13})$};
\node at (5,2.5) {$z_\eta\sigma(-A^{\frac13})$};
\node at (5,-2.5) {$\overline{z_\eta}\sigma(-A^{\frac13})$};
			\end{tikzpicture}
		\end{center}
\caption{Semi-lines contained the eigenvalues of $-\mathbb{B}_{(\eta)}$ and $0<\theta_\eta<\dfrac{\pi}{2}$, $ \eta >1$.}\label{figadadaf1New}
	\end{figure}

\subsection{Sectoriality}

Initially, we prove the following theorem on the sectoriality of the operator $\mathbb{B}_{(\eta)}$ for $\eta>1$.

\begin{theorem}\label{PosPDE}
Let $\eta>1$. The unbounded linear operator $\mathbb{B}_{(\eta)}$ defined in \eqref{Eq4aaabb}-\eqref{Eq4bbbbbbb} is  a sectorial operator.
\end{theorem}

\proof
In this proof, $M$ will denote a positive constant, not necessarily the same one. Let $\lambda\in\mathbb{C}$, then
\[
\lambda I-\mathbb{B}_{(\eta)}=\Big[\begin{smallmatrix} \lambda I- A^{\frac13} & I & 0 \\ 0 & \lambda I & I\\ 0 & - A^{\frac23} & \lambda I-(\eta-1) A^{\frac13} \end{smallmatrix}\Big]. 
\]
From Lemma\ref{MleNewddsd} it follows that 
\begin{equation}
\rho(\mathbb{B}_{(\eta)})=\rho(A^\frac13)\cap\rho(c_{\eta}A^\frac13)\cap\rho(d_{\eta}A^\frac13).
\end{equation}
Note that for $\lambda\in\rho(\mathbb{B}_{(\eta)})$ we have 
\begin{equation}\label{aw25}
(\lambda I-\mathbb{B}_{(\eta)})^{-1}=D_\eta(\lambda)^{-1}\left[\begin{smallmatrix}
(\lambda I-c_{\eta} A^{\frac13})(\lambda I-d_{\eta} A^{\frac13}) & -(\lambda I-(\eta-1)A^{\frac13})& I \\
0&  (\lambda I-A^{\frac13})(\lambda I-(\eta-1)A^{\frac13}) &  -(\lambda I-A^{\frac13})\\
0  & (\lambda I -A^{\frac13})A^{\frac23}& \lambda(\lambda I -A^{\frac13})
\end{smallmatrix}\right]
\end{equation}
where
\[
D_\eta(\lambda)=(\lambda I -A^{\frac13})(\lambda I -c_{\eta}A^{\frac13})(\lambda I -d_{\eta}A^{\frac13})
\]
with
\begin{eqnarray}
c_{\eta}&=&\frac12\left(\eta-1+\sqrt{\eta^2-2\eta-3}\right)\\
d_{\eta}&=&\frac12\left(\eta-1-\sqrt{\eta^2-2\eta-3}\right).
\end{eqnarray}

Since $Re(c_\eta)>0$ and $Re(d_\eta)>0$, for $\eta>1$, $c_\eta A^{\frac13}$ and $d_\eta A^{\frac13}$ are sectorial operators. Let $S_{A^\frac13}, S_{c_\eta A^{\frac13}}, S_{d_\eta A^\frac13} $ be sectors such that $S_{A^\frac13}\subset \rho(A^\frac13)$, $S_{c_\eta A^\frac13}\subset \rho(c_\eta A^\frac13)$, $S_{d_\eta A^\frac13}\subset \rho(d_\eta A^\frac13)$ and

\begin{eqnarray*}
\|\lambda I-A^\frac13\|_{\mathcal{L}(X)}&<&\frac{M}{|\lambda|}, \ \ \text{for each}\ \  \lambda \in S_{A^\frac13}\\
\|\lambda I-c_\eta A^\frac13\|_{\mathcal{L}(X)}&<&\frac{M}{|\lambda|}, \ \ \text{for each}\ \  \lambda \in S_{A^\frac13}\\
\|\lambda I-d_\eta A^\frac13\|_{\mathcal{L}(X)}&<&\frac{M}{|\lambda|}, \ \ \text{for each}\ \  \lambda \in S_{A^\frac13}
\end{eqnarray*}
for some $M>0$. We will prove that $\mathbb{B}_\eta$ is a sectorial operator using the sector 

$$S_{\mathbb{B}_{(\eta)}}:=S_{A^\frac13}\cap S_{c_\eta A^{\frac13}}\cap S_{d_\eta A^\frac13}.$$
 
It is immediate that $S_{\mathbb{B}_{(\eta)}}\subset \rho(\mathbb{B}_{(\eta)})$. If $\lambda\in S_{\mathbb{B}_{(\eta)}}$ and ${\bf u}=\left[\begin{smallmatrix}
u\\
v\\
w
\end{smallmatrix}\right]\in Z$ with $\|{\bf u}\|_Z \leqslant 1$, then writing 
\[
(\lambda I-\mathbb{B}_{(\eta)})^{-1}{\bf u}=\left[\begin{smallmatrix}
\varphi_1\\
\varphi_2\\
\varphi_3
\end{smallmatrix}\right]
\] 
where
\begin{equation}\label{masls9}
\begin{split}
\varphi_1&=(\lambda^2 I-(\eta-1)\lambda A^{\frac13}-A^{\frac23})  D_\eta(\lambda)^{-1}u + (-\lambda I+(\eta-1) A^{\frac13})  D_\eta(\lambda)^{-1}v + D_\eta(\lambda)^{-1}w,\\
\varphi_2&= (\lambda^2 I-\eta\lambda A^{\frac13}+(\eta-1)A^{\frac23}) D_\eta(\lambda)^{-1}v  +(-\lambda I+A^{\frac13}) D_\eta(\lambda)^{-1}w,\\
\varphi_3&= (\lambda A^{\frac23}-A) D_\eta(\lambda)^{-1}v + (\lambda^2 I-\lambda A^{\frac13})D_\eta(\lambda)^{-1}w.
\end{split}
\end{equation}

In order to conclude that 
$$\|\varphi_1\|_{X^{\frac23}}+\|\varphi_2\|_{X^{\frac13}}+\|\varphi_3\|_{X}< \frac{M}{|\lambda|}$$
We only need to show that  $\lambda A^\frac13 D_{\eta}(\lambda)^{-1}$, $ A^\frac23 D_{\eta}(\lambda)^{-1}$, $\lambda^2 D_{\eta}(\lambda)^{-1}$ $\in\mathcal{L}(X)$ and are bounded by $M/|\lambda|$, which is clear because $D_{\eta}(\lambda)^{-1}\in\mathcal{L}(X)$ and
$$\|D_{\eta}(\lambda)^{-1}\|_{\mathcal{L}(X)}<\frac{M}{|\lambda|^3}$$
\qed

\begin{lemma}
Let $Z_{-1}$ denote the extrapolation space of $Z$ generated by $\mathbb{B}_{(\eta)}$. Then
\[
Z_{-1}=X^{\frac13}\times X\times X^{-\frac13}.
\]
\end{lemma}

\proof
The extrapolation space of $Z$ is the completion of the normed space $(Z,\|\mathbb{B}_{(\eta)}^{-1}\|_Z)$. Since for $\left[\begin{smallmatrix}
u\\
v\\
w
\end{smallmatrix}\right]\in Z$ we have
\[
\mathbb{B}_{(\eta)}^{-1}\left[\begin{smallmatrix}
u\\
v\\
w
\end{smallmatrix}\right]=\Big[\begin{smallmatrix} \eta A^{-\frac13} u+ (1-\eta) A^{-\frac23} v+ A^{-1}w\\  (1-\eta) A^{-\frac13}v+ A^{-\frac23} w \\ v & \end{smallmatrix}\Big],
\]
and consequently
\[
\begin{split}
&\Big\|\mathbb{B}_{(\eta)}^{-1}\left[\begin{smallmatrix}
u\\
v\\
w
\end{smallmatrix}\right]\Big\|_{X^{\frac23}\times X^{\frac13}\times X}\\
&= \|\eta A^{-\frac13} u+ (1-\eta) A^{-\frac23} v+ A^{-1}w\|_{X^{\frac23}}+\|(1-\eta) A^{-\frac13}v+ A^{-\frac23} w \|_{X^{\frac13}}+\|v\|_X\\
&\leqslant \eta\| A^{\frac13} u\|_X+\max\{1,\eta-1\}\|v\|_X +2\|A^{-\frac{1}{3}}w\|_X\\
&\leqslant C_{1\eta}\Big\|\left[\begin{smallmatrix}
u\\
v\\
w
\end{smallmatrix}\right]\Big\|_{X\times X^{-\frac13}\times X^{-\frac23}},
\end{split}
\]
where $C_{1\eta}=\max\{2,\eta\}>0$.

By other hand
\[
\begin{split}
&\Big\|\left[\begin{smallmatrix}
u\\
v\\
w
\end{smallmatrix}\right]\Big\|_{X^{\frac{1}{3}}\times X\times X^{-\frac13}}\\
 &= \|u\|_{X^{\frac{1}{3}}}+\|v\|_X+\|w\|_{X^{-\frac{1}{3}}}\\
&= \dfrac{1}{\eta}\|\eta A^{-\frac{1}{3}}u\|_{X^{\frac{2}{3}}} +\|v\|_X+\|w\|_{X^{-\frac{1}{3}}}\\
&\leqslant \dfrac{1}{\eta}\|\eta A^{-\frac{1}{3}}u+ (1-\eta) A^{-\frac23} v+ A^{-1}w\|_{X^{\frac{2}{3}}}+\Big(2-\dfrac{1}{\eta}\Big)\|v\|_X+\Big(\dfrac{1}{\eta}+1\Big)\|w\|_{X^{-\frac{1}{3}}}\\
&\leqslant \dfrac{1}{\eta}\|\eta A^{-\frac{1}{3}}u+ (1-\eta) A^{-\frac23} v+ A^{-1}w\|_{X^{\frac{2}{3}}}+\Big(2+\eta-\dfrac{2}{\eta}\Big)\|v\|_X+\Big(\dfrac{1}{\eta}+1\Big)\|(1-\eta) A^{-\frac13} v+A^{-\frac{2}{3}}w\|_{X^{\frac{1}{3}}},\\
&\leqslant C_{2\eta}\Big\|\mathbb{B}_{(\eta)}^{-1}\left[\begin{smallmatrix}
u\\
v\\
w
\end{smallmatrix}\right]\Big\|_{X^{\frac23}\times X^{\frac13}\times X},
\end{split}
\]
for some $C_{2\eta}=2+\eta-\dfrac{2}{\eta}>0$.

Hence, there exist $C_{1\eta}>0$ and $C_{2\eta}>0$ such that 
\[
\Big\|\mathbb{B}_{(\eta)}^{-1}\left[\begin{smallmatrix}
u\\
v\\
w
\end{smallmatrix}\right]\Big\|_{X^{\frac23}\times X^{\frac13}\times X}\leqslant C_{1\eta} \Big\|\left[\begin{smallmatrix}
u\\
v\\
w
\end{smallmatrix}\right]\Big\|_{X\times X^{-\frac13}\times X^{-\frac23}}\leqslant C_{1\eta}C_{2\eta} \Big\|\mathbb{B}_{(\eta)}^{-1}\left[\begin{smallmatrix}
u\\
v\\
w
\end{smallmatrix}\right]\Big\|_{X^{\frac23}\times X^{\frac13}\times X},
\]
completions of $(X^{\frac23}\times X^{\frac13}\times X,\|\mathbb{B}_{(\eta)}^{-1}\cdot\|_{X^{\frac23}\times X^{\frac13}\times X})$ and $(X^{\frac23}\times X^{\frac13}\times X,\|\mathbb{B}_{(\eta)}^{-1}\cdot\|_{X^{\frac{1}{3}}\times X\times X^{-\frac13}})$ coincide.
\qed

Consider the closed extension of $\mathbb{B}_{(\eta)}$ to $Z_{-1}$ (see \cite[Page 262]{A}) and still denote it by $\mathbb{B}_{(\eta)}$. Then $\mathbb{B}_{(\eta)}$ is a sectorial and positive operator in $Z_{-1}$ with the domain $Z_{-1}^1=X^{\frac23}\times X^{\frac13}\times X$; the imaginary powers of $\mathbb{B}_{(\eta)}$ are bounded.  Our next concern will be to obtain embeddings of the spaces from the
fractional power scale $Z_{-1}^\alpha$, $\alpha\geqslant0$, generated by $(\mathbb{B}_{(\eta)},Z_{-1})$.

Before we can proceed we need the following general interpolation result.

\begin{proposition}
Let $\mathcal{X}_i$, $\mathcal{Y}_i$, $\mathcal{M}_i$, $i=1,2,3$ be the Banach spaces such that $\mathcal{X}_1\subset \mathcal{X}_0$, $\mathcal{Y}_1\subset \mathcal{Y}_0$, $\mathcal{M}_1\subset \mathcal{M}_0$ topologically and algebraically. Then
\[
[\mathcal{X}_0\times \mathcal{Y}_0\times \mathcal{M}_0,\mathcal{X}_1\times \mathcal{Y}_1\times \mathcal{M}_1]_\theta=[\mathcal{X}_0,\mathcal{X}_1]_\theta\times [\mathcal{Y}_0,\mathcal{Y}_1]_\theta\times [\mathcal{M}_0,\mathcal{M}_1]_\theta,
\]
for any $\theta\in(0,1)$.
\end{proposition}

\proof
The proof is an immediate consequence of the definition of complex interpolation spaces in \cite[Section 1.9.2]{HTr}.
\qed

The following result also may be established by a straightforward extension of   \cite[Theorem 2]{CCCC} so we omit its proof.

\begin{lemma}\label{asa9345y}
Denote by $\{Z_{-1}^\alpha: \alpha\in[0,1]\}$ the partial fractional power   scale generated by operator $\mathbb{B}_{(\eta)}$ in $Z_{-1}$. Then
\[
Z_{-1}^{k+\alpha}=X^{\frac{k+1+\alpha}{3}}\times X^{\frac{k+\alpha}{3}}\times X^{\frac{k-1+\alpha}{3}},\quad \alpha\in[0,1],\ k\in\mathbb{N}.
\]
\end{lemma}

For better understanding the relation of the fractional power scale-spaces of the operator $\mathbb{B}_{(\eta)}$, we construct the following diagram.
\begin{figure}[h]\label{Figure_A0100N}
\begin{tikzpicture}
\draw[->] (-3.5,0) -- (3.5,0);
\node at (0,-0.5) {$Z$};
\filldraw[gray!40!,fill opacity=0.2] (-2.2,-0.3) -- (-1.8,-0.1) -- (-1.8,0.4) -- (-2.2,0.2) -- cycle;
\filldraw[gray!40!,fill opacity=0.2] (-0.2,-0.3) -- (0.2,-0.1) -- (0.2,0.4) -- (-0.2,0.2) -- cycle;
\filldraw[gray!40!,fill opacity=0.2] (1.8,-0.3) -- (2.2,-0.1) -- (2.2,0.4) -- (1.8,0.2) -- cycle;
\filldraw[gray!40!,fill opacity=0.2] (-2.2,-1.8) -- (-1.8,-1.6) -- (-1.8,-1.1) -- (-2.2,-1.3) -- cycle;
\filldraw[gray!40!,fill opacity=0.2] (-0.2,-1.8) -- (0.2,-1.6) -- (0.2,-1.1) -- (-0.2,-1.3) -- cycle;
\filldraw[gray!40!,fill opacity=0.2] (1.8,-1.8) -- (2.2,-1.6) -- (2.2,-1.1) -- (1.8,-1.3) -- cycle;
\node at (-2,-0.5) {$Z^1$};
\node at (-3,-0.3) {$\cdots$};
\node at (2,-0.5) {$(Z^1)'=Z^{-1}$};
\node at (4,-0.3) {$\cdots$};
\node at (-1,0.6) {$\mathbb{B}_{(\eta)}$};
\draw[densely dotted,-latex'] (-2,0.1) to [out=30,in=150] (0,0.1);
\draw[->] (-3.5,-1.5) -- (3.5,-1.5);
\node at (0,-2.1) {$Z^1_{-1}$};
\node at (-2,-2.1) {$Z^\alpha_{-1}$};
\node at (-3,-1.85) {$\cdots$};
\node at (2,-2.1) {$Z_{-1}$};
\node at (4,-1.85) {$\cdots$};
\node at (1,-0.95) {$\mathbb{B}_{(\eta)}$};
\draw[densely dotted,-latex'] (0,-1.45) to [out=30,in=150] (2,-1.45);
\end{tikzpicture}
\caption{Fractional power scale generated by operator $\mathbb{B}_{(\eta)}$ in $Z_{-1}$ and $\alpha>1$.}
\end{figure}
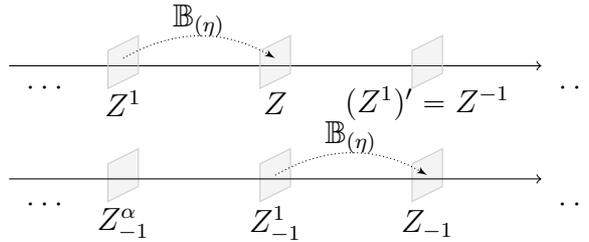

In this section, we consider the case $\eta>1$. Now we prove some of the main results of this paper.

\begin{proposition}
Let $Z_{-1}$ denote the extrapolation space of $Z$ generated by $\mathbb{A}_\eta$. Then
\[
Z_{-1}=X^{\frac13}\times X\times X^{-\frac13}.
\]
\end{proposition}

\begin{figure}[h]\label{Figure_B}
\begin{tikzpicture}
\draw[->] (-3.5,0) -- (3.5,0);
\filldraw[gray!40!,fill opacity=0.2] (-2.2,-0.3) -- (-1.8,-0.1) -- (-1.8,0.4) -- (-2.2,0.2) -- cycle;
\filldraw[gray!40!,fill opacity=0.2] (-0.2,-0.3) -- (0.2,-0.1) -- (0.2,0.4) -- (-0.2,0.2) -- cycle;
\filldraw[gray!40!,fill opacity=0.2] (1.8,-0.3) -- (2.2,-0.1) -- (2.2,0.4) -- (1.8,0.2) -- cycle;
\filldraw[gray!40!,fill opacity=0.2] (-2.2,-1.8) -- (-1.8,-1.6) -- (-1.8,-1.1) -- (-2.2,-1.3) -- cycle;
\filldraw[gray!40!,fill opacity=0.2] (-0.2,-1.8) -- (0.2,-1.6) -- (0.2,-1.1) -- (-0.2,-1.3) -- cycle;
\filldraw[gray!40!,fill opacity=0.2] (1.8,-1.8) -- (2.2,-1.6) -- (2.2,-1.1) -- (1.8,-1.3) -- cycle;
\node at (0,-0.5) {$Z$};
\node at (-2,-0.5) {$Z^1$};
\node at (2,-0.5) {$(Z^1)'=Z^{-1}$};
\node at (-1,0.65) {$\mathbb{A}_\eta$};
\draw[densely dotted,-latex'] (-2,0.1) to [out=30,in=150] (0,0.1);
\draw[->] (-3.5,-1.5) -- (3.5,-1.5);
\node at (0,-2.1) {$Z^1_{-1}$};
\node at (-2,-2.1) {$Z^\alpha_{-1}$};
\node at (2,-2.1) {$Z_{-1}$};
\node at (1,-0.9) {$\mathbb{A}_\eta$};
\draw[densely dotted,-latex'] (0,-1.45) to [out=30,in=150] (2,-1.45);
\node at (-3,-0.3) {$\cdots$};
\node at (-3,-1.85) {$\cdots$};
\node at (4,-1.85) {$\cdots$};
\node at (4,-0.3) {$\cdots$};
\end{tikzpicture}
\caption{Fractional power scale generated by operator $\mathbb{A}_\eta$ in $Z_{-1}$ and $\alpha>1$.}
\end{figure}
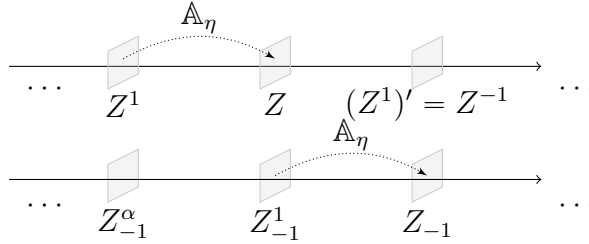

\proof
Again, following \cite{A} and \cite{CC_2} we recall that the extrapolation space of $Z$ is the completion of the normed space $(Z,\|\mathbb{A}_\eta^{-1}\|_Z)$. Since for $\Big[\begin{smallmatrix} u\\ v\\ w \end{smallmatrix}\Big]\in Z$
\[
\Big\|\mathbb{A}_\eta^{-1}\Big[\begin{smallmatrix} u\\ v\\ w \end{smallmatrix}\Big]\Big\|_{X^{\frac23}\times X^{\frac13}\times X}=\Big\|\Big[\begin{smallmatrix}  \eta A^{-\frac13}u+\eta A^{-\frac23} v+A^{-1}w\\ -u \\ -v \end{smallmatrix}\Big]\Big\|_{X^{\frac23}\times X^{\frac13}\times X}\leqslant 3(1+\eta) \Big\|\Big[\begin{smallmatrix} u\\ v\\ w \end{smallmatrix}\Big]\Big\|_{X^{\frac13}\times X\times X^{-\frac13}},
\]
and
\[
 \Big\|\Big[\begin{smallmatrix} u\\ v\\ w \end{smallmatrix}\Big]\Big\|_{X^{\frac13}\times X\times X^{-\frac13}}\leqslant 3(1+\eta) \Big\|\mathbb{A}_\eta^{-1}\Big[\begin{smallmatrix} u\\ v\\ w \end{smallmatrix}\Big]\Big\|_{X^{\frac23}\times X^{\frac13}\times X},
\]
completions of $(X^{\frac23}\times X^{\frac13}\times X,\|\mathbb{A}_\eta^{-1}\cdot\|_{X^{\frac23}\times X^{\frac13}\times X})$ and $(X^{\frac23}\times X^{\frac13}\times X,\|\cdot\|_{X^{\frac13}\times X\times X^{-\frac13}})$ coincide, see the Figure \ref{Figure_B}.
\qed

Consider the closed extension of $\mathbb{A}_\eta$ to $Z_{-1}$ (see \cite[Page 262]{A}) and still denote it by $\mathbb{A}_\eta$  with the domain $Z_{-1}^1=X^{\frac23}\times X^{\frac13}\times X$. 

For better understanding the relation of the fractional power scale-spaces of the operator $\mathbb{A}_\eta$, we construct the following diagram.

\subsection{Analysis by reducing the order}

Namely, we consider the reduction of order 
\begin{equation}\label{cv1}
v=u_t + A^{\frac13}u
\end{equation}
for positive time, where $u$  is an unknown function to be determined, and we obtain the following equation of second order 
\begin{equation}\label{Eq2}
v_{tt} + A^{\frac23} v+ (\eta-1)A^{\frac13}v_t  =g(v)
\end{equation}
for positive time, where $\eta>0$ and 
\begin{equation}\label{Eq2g}
g(v)=f(u).
\end{equation}
 
Note that we can reviews \eqref{Eq2} as follows 
\begin{equation}\label{E00q1Linear31}
\left[ \begin{matrix} v \\ v_t \end{matrix} \right]_t + \left[ \begin{matrix} 0 & -I \\ A^\frac23 & (\eta-1)A^\frac13  \end{matrix} \right]\left[ \begin{matrix} v \\ v_t \end{matrix} \right] =\left[ \begin{matrix} 0 \\ g(v) \end{matrix} \right], \quad \left[ \begin{matrix} v(0) \\ v_t(0) \end{matrix} \right]=\left[ \begin{matrix} v_0 \\ w_0 \end{matrix} \right]\in X^\frac13 \times X.
\end{equation}

Solving \eqref{E00q1Linear31} we find $e^{-\Lambda t}\left[ \begin{matrix} v_0 \\ w_0 \end{matrix} \right]=\left[ \begin{matrix} v(t,v_0,w_0) \\ v_t (t,v_0,w_0)\end{matrix} \right]\in X^\frac13 \times X$
where $\Lambda:D(\Lambda)\subset X^\frac13\times X\to X^\frac13\times X$ is defined by,  $D(\Lambda)=X^\frac23\times X^\frac13$
$$
\Lambda \left[ \begin{matrix} v_0 \\ w_0 \end{matrix} \right]= \left[ \begin{matrix} 0 & -I \\ A^\frac23 & (\eta-1)A^\frac13 \end{matrix} \right]\left[ \begin{matrix} v_0 \\ w_0 \end{matrix} \right]=\left[ \begin{matrix} -w_0 \\ A^\frac23 v_0 \end{matrix} \right]
$$

It follows that $\mathbb{R}^+\ni t\mapsto  v(t,v_0,w_0)\in X^\frac13$ is a continuous function. We can then try to solve
\begin{equation}\label{cv1_new}
u_t + A^{\frac13}u=v(t,v_0,w_0), \quad u(0)=u_0
\end{equation}
through the variation of constants formula
\begin{equation}\label{cv21}
u(t)=e^{-A^{\frac13}t}u(0)+\displaystyle\int_0^t e^{-A^{\frac13}(t-s)}v(s)ds.
\end{equation}
for positive time.

We next solve \eqref{Eq2} following the same ideas in \cite{CC_2,CC_3,CT1} and \cite{ChD}. Next, we use \eqref{cv21} to obtain a local solution to our original differential equation. A similar argument has been used in \cite{BP} and \cite{KLM} for the classical MGT equation. 

To better present our results we introduce some notations. We consider the space $Y=X^{\frac13}\times X$ equipped the norm given by
\[
\Big\|\begin{bmatrix} v\\ w \end{bmatrix}\Big\|^2_{Y}=\|v\|^2_{X^{\frac13}}+\|w\|^2_X,\ \forall \begin{bmatrix} v\\ w \end{bmatrix}\in Y.
\] 
We can rewrite the initial value problem associated with equation \eqref{Eq2}  as the Cauchy problem in $Y$
\begin{equation}\label{Eq3a}
\dfrac{d}{dt}\begin{bmatrix} v\\ w \end{bmatrix}+ \varLambda\begin{bmatrix} v\\ w \end{bmatrix}=G(\begin{bmatrix} v\\ w \end{bmatrix}),\ t>0,
\end{equation}
and
\begin{equation}\label{Eq3b}
\begin{bmatrix} v\\ w \end{bmatrix}(0)=\begin{bmatrix} v_0\\ w_0 \end{bmatrix}
\end{equation}
where $w=v_t$ and the unbounded linear operator $\varLambda:D(\varLambda)\subset Y\to Y$ is defined by
\[
D(\varLambda)=X^{\frac23}\times X^{\frac13}
\]
and
\[
\varLambda \begin{bmatrix} v\\ w \end{bmatrix}:=\begin{bmatrix} 0 & -I\\  A^{\frac23} & (\eta-1)A^{\frac13} \end{bmatrix}\begin{bmatrix} v\\ w \end{bmatrix}=\begin{bmatrix} -w\\ A^{\frac23}v+(\eta-1)A^{\frac13}w \end{bmatrix},\ \forall\begin{bmatrix} v\\ w \end{bmatrix}\in X^{\frac23}\times X^{\frac13}.
\]

From now on, we denote 
\[
Y^1=D(\varLambda)=X^{\frac23}\times X^{\frac13}.
\]

The nonlinearity $G$ is given by 
\begin{equation}\label{DefinitionG}
G(\begin{bmatrix} v\\ w \end{bmatrix})=\begin{bmatrix} 0\\ g(v) \end{bmatrix}
\end{equation}
where $g$ is given by \eqref{Eq2g}.

The following result may be established by a straightforward extension of \cite[Lemma 1]{CC_2} so we omit its proof.

\begin{lemma}\label{LL}
The following conditions hold.
\begin{itemize}
\item[$i)$] The linear operator $\varLambda$ is closed and densely defined; 
\item[$ii)$] $0\in\rho(\varLambda)$ with $\varLambda^{-1}=\Big[\begin{smallmatrix} (\eta-1)A^{-\frac13} & A^{-\frac23}\\ -I  & 0 \end{smallmatrix}\Big]$. Moreover,  $\varLambda$ has compact resolvent;
\item[$iii)$] The spectrum of $\varLambda$ consists of isolated eigenvalues $\lambda^\pm_n$ given by
\[
\lambda^\pm_n=(\eta\pm\sqrt{\eta^2-1})\sqrt[3]{\mu_n},\ \mu_n\in\sigma(A);
\] 
\item[$iv)$] The linear operator $-\varLambda$ generates in $Y$ a $C^0$ analytic semigroup $\{e^{-\varLambda t}:t\geqslant0\}$;
\item[$v)$] The semigroup $\{e^{-\varLambda t}:t\geqslant0\}$ in $Y$ are compact and asymptotically decaying;
\item[$vi)$] For each $\begin{bmatrix} v\\ w \end{bmatrix}\in Y^1$, we have
\[
\dfrac{1}{\eta}\Big\|\begin{bmatrix} v\\ w \end{bmatrix}\Big\|_{Y^1}\leqslant \Big\|\varLambda\begin{bmatrix} v\\ w \end{bmatrix}\Big\|_{Y}\leqslant \eta\Big\|\begin{bmatrix} v\\ w \end{bmatrix}\Big\|_{Y^1}.
\]
\end{itemize}
\end{lemma}

The Lemma \ref{LL} allows us to define the fractional power $\varLambda^{-\alpha}$ of order $\alpha\in(0,1)$. Denote by $Y^\alpha=D(\varLambda^{\alpha})$ for $\alpha\in[0,1)$, taking $\varLambda^0:=I$ on $Y^0:=Y$ when $\alpha=0$. Recall that $Y^{\alpha}$ is dense in $Y$  for all $\alpha\in(0,1]$, for details see  \cite[Theorem 4.6.5]{A}. The fractional power space $X^\alpha$ endowed with the norm 
\[
\|\cdot\|_{Y^\alpha}:=\|\varLambda^{\alpha} \cdot\|_Y
\] 
is a Banach space. 

\begin{remark}\label{expspace}
Let $X_{-1}$ denote the extrapolation space of $X$ generated by $A^{\frac13}$. Consider the closed extension of $A^{\frac13}$ to $X_{-1}$ (see \cite[Page 262]{A}) and still denote it by $A^{\frac13}$. Then $A^{\frac13}$ is a sectorial and positive operator in $X_{-1}$ with the domain $X_{-1}^{\frac13}$; the imaginary powers of $A^{\frac13}$ are bounded, see e.g. \cite[Proposition 5]{CCCC}.  Our next concern will be to obtain embeddings of the spaces from the
fractional power scale $X_{-1}^{\frac{\alpha}{3}}$, $\alpha\geqslant0$, generated by $(A^{\frac13},X_{-1})$.
\end{remark}

For better understanding the relation of the fractional power scale-spaces of the operator $A^{\frac13}$, we construct the following diagram.
\begin{figure}[h]
\begin{tikzpicture}
\draw[->] (-3.5,0) -- (3.5,0);
\filldraw[gray!40!,fill opacity=0.2] (-2.2,-0.3) -- (-1.8,-0.1) -- (-1.8,0.4) -- (-2.2,0.2) -- cycle;
\filldraw[gray!40!,fill opacity=0.2] (-0.2,-0.3) -- (0.2,-0.1) -- (0.2,0.4) -- (-0.2,0.2) -- cycle;
\filldraw[gray!40!,fill opacity=0.2] (1.8,-0.3) -- (2.2,-0.1) -- (2.2,0.4) -- (1.8,0.2) -- cycle;
\filldraw[gray!40!,fill opacity=0.2] (-2.2,-1.8) -- (-1.8,-1.6) -- (-1.8,-1.1) -- (-2.2,-1.3) -- cycle;
\filldraw[gray!40!,fill opacity=0.2] (-0.2,-1.8) -- (0.2,-1.6) -- (0.2,-1.1) -- (-0.2,-1.3) -- cycle;
\filldraw[gray!40!,fill opacity=0.2] (1.8,-1.8) -- (2.2,-1.6) -- (2.2,-1.1) -- (1.8,-1.3) -- cycle;
\node at (0,-0.5) {$X$};
\node at (-2,-0.5) {$X^{\frac13}$};
\node at (2,-0.5) {$(X^{\frac13})'=X^{-\frac13}$};
\node at (-1,0.65) {$A^{\frac13}$};
\draw[densely dotted,-latex'] (-2,0.1) to [out=30,in=150] (0,0.1);
\draw[->] (-3.5,-1.5) -- (3.5,-1.5);
\node at (0,-2.1) {$X^{\frac13}_{-1}$};
\node at (-2,-2.1) {$X^{\frac{\alpha}{3}}_{-1}$};
\node at (2,-2.1) {$X_{-1}$};
\node at (1,-0.9) {$A^{\frac13}$};
\draw[densely dotted,-latex'] (0,-1.45) to [out=30,in=150] (2,-1.45);
\node at (-3,-0.3) {$\cdots$};
\node at (-3,-1.85) {$\cdots$};
\node at (4,-1.85) {$\cdots$};
\node at (4,-0.3) {$\cdots$};
\end{tikzpicture}
\caption{Fractional power scale generated by operator $A^{\frac13}$ in $X_{-1}$ and $\alpha>1$.}
\end{figure}

The following result also may be established by a straightforward extension of   \cite[Proposition 3]{CCCC} or \cite[Theorem 2.3]{CC_3} so we omit its proof.

\begin{theorem}\label{TheoCCD}
For each $\alpha\in[0,1]$, the fractional power spaces $Y^\alpha$ associated to the operator $\varLambda$ coincide with $X^{\frac{1+\alpha}{3}}\times X^{\frac{\alpha}{3}}$  with equivalent norms.
\end{theorem}

\begin{lemma}
Let $Y_{-1}$ denote the extrapolation space of $Y$ generated by $\varLambda$. Then
\[
Y_{-1}=X\times X^{-\frac13}.
\]
\end{lemma}

\proof
Following \cite{A} and \cite{CC_2} we recall that the extrapolation space of $Y$ is the completion of the normed space $(Y,\|\varLambda^{-1}\|_Y)$. Since for $\begin{bmatrix} v\\ w \end{bmatrix}\in Y$
\[
\Big\|\varLambda^{-1}\begin{bmatrix} v\\ w \end{bmatrix}\Big\|_{X^{\frac13}\times X}\leqslant\eta \Big\|\begin{bmatrix} v\\ w \end{bmatrix}\Big\|_{X\times X^{-\frac13}}\leqslant\eta^2 \Big\|\varLambda^{-1}\begin{bmatrix} v\\ w \end{bmatrix}\Big\|_{X^{\frac13}\times X},
\]
completions of $(X^{\frac13}\times X,\|\varLambda^{-1}\cdot\|_{X^{\frac13}\times X})$ and $(X^{\frac13}\times X,\|\varLambda^{-1}\cdot\|_{X\times X^{-\frac13}})$ coincide, see the Figure \ref{Figure_A01}.
\qed

Consider the closed extension of $\varLambda$ to $Y_{-1}$ (see \cite[Page 262]{A}) and still denote it by $\varLambda$. Then $\varLambda$ is a sectorial and positive operator in $Y_{-1}$ with the domain $Y_{-1}^1=X^{\frac13}\times X$; the imaginary powers of $\varLambda$ are bounded, see e.g. \cite[Proposition 5]{CCCC}.  Our next concern will be to obtain embeddings of the spaces from the
fractional power scale $Y_{-1}^\alpha$, $\alpha\geqslant0$, generated by $(\varLambda,Y_{-1})$.

The following result also may be established by a straightforward extension of   \cite[Theorem 2]{CCCC} so we omit its proof.

\begin{lemma}\label{asa9345y}
Denote by $\{Y_{-1}^\alpha: \alpha\in[0,1]\}$ the partial fractional power   scale generated by operator $\varLambda$ in $Y_{-1}$. Then
\[
Y_{-1}^{k+\alpha}=X^{\frac{k+\alpha}{3}}\times X^{\frac{k-1+\alpha}{3}},\quad \alpha\in[0,1],\ k\in\mathbb{N}.
\]
\end{lemma}

For better understanding the relation of the fractional power scale-spaces of the operator $\varLambda$, we construct the following diagram.
\begin{figure}[h]\label{Figure_A01}
\begin{tikzpicture}
\draw[->] (-3.5,0) -- (3.5,0);
\filldraw[gray!40!,fill opacity=0.2] (-2.2,-0.3) -- (-1.8,-0.1) -- (-1.8,0.4) -- (-2.2,0.2) -- cycle;
\filldraw[gray!40!,fill opacity=0.2] (-0.2,-0.3) -- (0.2,-0.1) -- (0.2,0.4) -- (-0.2,0.2) -- cycle;
\filldraw[gray!40!,fill opacity=0.2] (1.8,-0.3) -- (2.2,-0.1) -- (2.2,0.4) -- (1.8,0.2) -- cycle;
\filldraw[gray!40!,fill opacity=0.2] (-2.2,-1.8) -- (-1.8,-1.6) -- (-1.8,-1.1) -- (-2.2,-1.3) -- cycle;
\filldraw[gray!40!,fill opacity=0.2] (-0.2,-1.8) -- (0.2,-1.6) -- (0.2,-1.1) -- (-0.2,-1.3) -- cycle;
\filldraw[gray!40!,fill opacity=0.2] (1.8,-1.8) -- (2.2,-1.6) -- (2.2,-1.1) -- (1.8,-1.3) -- cycle;
\node at (0,-0.5) {$Y$};
\node at (-2,-0.5) {$Y^1$};
\node at (-3,-0.3) {$\cdots$};
\node at (2,-0.5) {$(Y^1)'=Y^{-1}$};
\node at (4,-0.3) {$\cdots$};
\node at (-1,0.6) {$\varLambda$};
\draw[densely dotted,-latex'] (-2,0.1) to [out=30,in=150] (0,0.1);
\draw[->] (-3.5,-1.5) -- (3.5,-1.5);
\node at (0,-2.1) {$Y^1_{-1}$};
\node at (-2,-2.1) {$Y^\alpha_{-1}$};
\node at (-3,-1.85) {$\cdots$};
\node at (2,-2.1) {$Y_{-1}$};
\node at (4,-1.85) {$\cdots$};
\node at (1,-0.95) {$\varLambda$};
\draw[densely dotted,-latex'] (0,-1.45) to [out=30,in=150] (2,-1.45);
\end{tikzpicture}
\caption{Fractional power scale generated by operator $\varLambda$ in $Y_{-1}$ and $\alpha>1$.}
\end{figure}

Since  the semigroup $\{e^{-\varLambda t}:t\geqslant0\}$ generated by $-\varLambda$ in $Y_{-1}$ is analytic and the linear Cauchy problem \ref{Eq3a}-\ref{Eq3b} with $\begin{bmatrix} v_0\\ w_0 \end{bmatrix}\in Y_{-1}$ has a unique solution
\begin{equation}\label{Sem_1}
\begin{bmatrix} v\\ w \end{bmatrix}(t)=e^{-\varLambda t}\begin{bmatrix} v_0\\ w_0 \end{bmatrix},\ t\geqslant0.
\end{equation}

In the following theorem we explain the smoothing action of the solution to the linear Cauchy problem \ref{Eq3a}-\ref{Eq3b}, it follows from \cite[Theorem 2.4]{CC_3} and \cite[Theorem 2.1.1]{ChD}. 

\begin{theorem}\label{TheoCCD2}
If $t>0$ and $\begin{bmatrix} v_0\\ w_0 \end{bmatrix}\in Y_{-1}$, then 
\[
\begin{bmatrix} v\\ w \end{bmatrix}(t)=e^{-\varLambda t}\begin{bmatrix} v_0\\ w_0 \end{bmatrix}\in Y_{-1}^\alpha,\ \mbox{for each}\ \alpha\geqslant0.
\]
In particular, if $t>0$ and $\begin{bmatrix} v_0\\ w_0 \end{bmatrix}\in Y^1_{-1}$, then 
\[
\begin{bmatrix} v\\ w \end{bmatrix}(t)=e^{-\varLambda t}\begin{bmatrix} v_0\\ w_0 \end{bmatrix}\in Y_{-1}^1.
\]
Moreover, 
\[
\begin{bmatrix} v\\ w \end{bmatrix}(\cdot)\in C^1((0,\infty), Y^\alpha_{-1}),\ \mbox{for each}\ \alpha\in [0,1).
\]
In particular, if  $\begin{bmatrix} v_0\\ w_0 \end{bmatrix}\in Y^1_{-1}$, then 
\[
\begin{bmatrix} v\\ w \end{bmatrix}(\cdot)\in C^1([0,\infty), Y_{-1}).
\]
\end{theorem}

\begin{theorem}\label{theo_fg}
Assume that $f$ is $\epsilon-$regular map relative to the pair $(X^{\frac13}, X)$ for $\epsilon\geqslant0$, in the sense of \eqref{asaE4}, then $g$ is $\epsilon-$regular map relative to the pair $(X_{-1}^\frac13,X_{-1})$, where $\{X_{-1}^{\frac{\alpha}{3}}:\alpha\geqslant0\}$ denotes the fractional power scale generated by operator $A^{\frac13}$ in $X_{-1}$.
\end{theorem}

\proof
We will prove that there exist  constants $c>0$, $\rho>1$, $\gamma(\epsilon)>0$ with $\rho\epsilon\leqslant \gamma(\epsilon)<\frac13$ such that $g:X_{-1}^{\frac13+\epsilon}\to X_{-1}^{\gamma(\epsilon)}$ and  
\[
\|g(\phi_1)-g(\phi_2)\|_{X_{-1}^{\gamma(\epsilon)}}\leqslant c\|\phi_1-\phi_2\|_{X_{-1}^{\frac13+\epsilon}}(1+\|\phi_1\|_{X_{-1}^{\frac13+\epsilon}}^{\rho-1}+\|\phi_2\|_{X_{-1}^{\frac13+\epsilon}}^{\rho-1}),
\]
for any $\phi_1,\phi_2\in X_{-1}^{\frac13+\epsilon}$, see the figures below.

For a better understanding of the $\epsilon$-regular map relative to the pair $(X_{-1}^{\frac13},X_{-1})$ for $\epsilon\geqslant0$, we construct the following diagram.
\begin{figure}[h]
\begin{tikzpicture}
\draw[->] (-3.5,0) -- (3,0);
\node at (1,0.01) {$\nshortmid$};
\node at (1,-0.3) {$X_{-1}$};
\node at (-2,0.01) {$\nshortmid$};
\node at (-2,-0.3) {$X_{-1}^\frac13$};
\node at (-0.5,0.8) {$g$};
\draw[densely dotted,-latex'] (-2,0.1) to [out=30,in=150] (1,0.1);
\draw[->] (-3.5,-1.5) -- (3,-1.5);
\node at (0.5,-1.51) {$\nshortmid$};
\node at (0.5,-1.85) {$X_{-1}^{\gamma(\epsilon)}$};
\node at (-2.7,-1.51) {$\nshortmid$};
\node at (-2.5,-1.85) {$X_{-1}^{\frac{1}{3}+\epsilon}$};
\node at (-1,-0.65) {$g$};
\draw[densely dotted,-latex'] (-2.7,-1.51) to [out=30,in=150] (0.5,-1.45);
\node at (-3,-0.3) {$\cdots$};
\node at (-3.5,-1.85) {$\cdots$};
\node at (1.5,-1.85) {$\cdots$};
\node at (2,-0.3) {$\cdots$};
\end{tikzpicture}
\caption{$X_{-1}^{\frac{1}{3}+\epsilon}\subset X_{-1}^{\frac13}$ and $X_{-1}^{\gamma(\epsilon)}\subset X_{-1}$}
\end{figure}
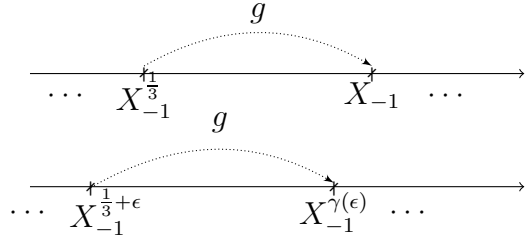

Since $f$ is an $\epsilon$-regular map relative to the pair $(X^{\frac13},X)$ for $\epsilon\geqslant0$, we have the following: if $g(\phi_i)=f(\psi_i)$, where 
\[
\psi_i=e^{-A^{\frac13}t}u(0)+\displaystyle\int_0^t e^{-A^{\frac13}(t-s)}\phi_i(s)ds,\quad i=1,2,
\]
and $\{e^{-A^{\frac13}t};t\geqslant0\}$ denotes the  $C^0-$semigroup corresponding to \eqref{cv1} with $v\equiv0$ subject to initial condition $u(0)=u_0\in X^{\frac{1}{3}+\epsilon}$ then for each $t\geqslant0$, $g=f\circ e^{-A^{\frac13}t}$ is an $\epsilon$-regular map relative to the pair $(X_{-1}^{\frac13},X_{-1})$ for $\epsilon\geqslant0$ because by Remark \ref{expspace} we have
\[
\begin{split}
&\|g(\phi_1)-g(\phi_2)\|_{X_{-1}^{\gamma(\epsilon)}}\\
&\leqslant c_0\|f(\psi_1)-f(\psi_2)\|_{X^{\gamma(\epsilon)}}\\
&\leqslant c_0\|\psi_1-\psi_2\|_{X^{\frac13+\epsilon}}(1+\|\psi_1\|_{X_{-1}^{\frac13+\epsilon}}^{\rho-1}+\|\psi_2\|_{X^{\frac13+\epsilon}}^{\rho-1})\\ 
&=c_1\Big\| \displaystyle\int_0^t e^{-A^{\frac13}(t-s)}\phi_1(s)ds-\displaystyle\int_0^t e^{-A^{\frac13}(t-s)}\phi_2(s)ds \Big\|_{X^{\frac13+\epsilon}}\times\\
&\times \Big(    1+\Big\|e^{-A^{\frac13}t}u(0)+\displaystyle\int_0^t e^{-A^{\frac13}(t-s)}\phi_1(s)ds\Big\|_{X_{-1}^{\frac13+\epsilon}}^{\rho-1}+\Big\|e^{-A^{\frac13}t}u(0)+\displaystyle\int_0^t e^{-A^{\frac13}(t-s)}\phi_2(s)ds\Big\|_{X^{\frac13+\epsilon}}^{\rho-1} \Big)\\
&=c_1\Big\| \displaystyle\int_0^t e^{-A^{\frac13}(t-s)}(\phi_1-\phi_2)(s)ds \Big\|_{X^{\frac13+\epsilon}}\times\\
&\times \Big(    1+\Big\|e^{-A^{\frac13}t}u(0)+\displaystyle\int_0^t e^{-A^{\frac13}(t-s)}\phi_1(s)ds\Big\|_{X_{-1}^{\frac13+\epsilon}}^{\rho-1}+\Big\|e^{-A^{\frac13}t}u(0)+\displaystyle\int_0^t e^{-A^{\frac13}(t-s)}\phi_2(s)ds\Big\|_{X^{\frac13+\epsilon}}^{\rho-1} \Big)\\
&\leqslant c_2\|\phi_1-\phi_2\|_{X_{-1}^{\frac13+\epsilon}}(1+\|\phi_1\|_{X_{-1}^{\frac13+\epsilon}}^{\rho-1}+\|\phi_2\|_{X_{-1}^{\frac13+\epsilon}}^{\rho-1})
\end{split}
\]
for some positive constants $c_0,c_1$ and $c_2$.
\qed

Thanks to \cite[Corollary 2]{CCCC} we have the following on the $\epsilon-$regular map $G$ defined in \eqref{DefinitionG}.

\begin{proposition}
For $\epsilon\geqslant0$, $G$ is $\epsilon-$regular map relative to the pair $(Y^1_{-1},Y_{-1})$; that is,  there exist  constants $C>0$, $\rho>1$, $\gamma(\epsilon)$ with $\rho\epsilon\leqslant \gamma(\epsilon)<1$ such that such that $G:Y_{-1}^{1+\epsilon}\to Y_{-1}^{\gamma(\epsilon)}$ and  
\begin{equation}\label{asa034mt56}
\Big\|G(\begin{bmatrix} \phi\\ \varphi \end{bmatrix})-G(\begin{bmatrix} \phi'\\ \varphi' \end{bmatrix})\Big\|_{Y_{-1}^{\gamma(\epsilon)}}\leqslant C\Big\|\begin{bmatrix} \phi\\ \varphi \end{bmatrix}-\begin{bmatrix} \phi'\\ \varphi' \end{bmatrix}\Big\|_{Y_{-1}^{1+\epsilon}}\Big(1+\Big\|\begin{bmatrix} \phi\\ \varphi \end{bmatrix}\Big\|_{Y_{-1}^{1+\epsilon}}^{\rho-1}+\Big\|\begin{bmatrix} \phi'\\ \varphi' \end{bmatrix}\Big\|_{Y_{-1}^{1+\epsilon}}^{\rho-1}\Big),
\end{equation}
for any $\begin{bmatrix} \phi\\ \varphi \end{bmatrix},\begin{bmatrix} \phi'\\ \varphi' \end{bmatrix}\in Y_{-1}^{1+\epsilon}.$
\end{proposition}

\proof
Let $\begin{bmatrix} \phi\\ \varphi \end{bmatrix},\begin{bmatrix} \phi'\\ \varphi' \end{bmatrix}\in Y_{-1}^{1+\epsilon}$. Thanks to Lemma \ref{asa9345y} we have 
\[
\Big\|G(\begin{bmatrix} \phi\\ \varphi \end{bmatrix})-G(\begin{bmatrix} \phi'\\ \varphi' \end{bmatrix})\Big\|_{Y_{-1}^{\gamma(\epsilon)}}=\Big\|\begin{bmatrix} 0\\ g(\phi)-g(\phi') \end{bmatrix}\Big\|_{Y_{-1}^{\gamma(\epsilon)}}=\|g(\phi)-g(\phi')\|_{X^{\frac{-1+\gamma(\epsilon)}{3}}}
\]
and from fractional power scale generated by operator $A^{\frac13}$ in $X_{-1}$ we obtain
\[
\Big\|G(\begin{bmatrix} \phi\\ \varphi \end{bmatrix})-G(\begin{bmatrix} \phi'\\ \varphi' \end{bmatrix})\Big\|_{Y_{-1}^{\gamma(\epsilon)}}\leqslant c\|g(\phi)-g(\phi')\|_{X_{-1}^{\frac{\gamma(\epsilon)}{3}}}
\]
for some $c>0$, and by Theorem \ref{theo_fg} we get
\[
\begin{split}
\Big\|G(\begin{bmatrix} \phi\\ \varphi \end{bmatrix})-G(\begin{bmatrix} \phi'\\ \varphi' \end{bmatrix})\Big\|_{Y_{-1}^{\gamma(\epsilon)}}&\leqslant  C\|\phi-\phi'\|_{X_{-1}^{\frac13+\tilde\epsilon}}(1+\|\phi\|_{X_{-1}^{\frac13+\tilde\epsilon}}^{\rho-1}+\|\phi'\|_{X_{-1}^{\frac13+\tilde\epsilon}}^{\rho-1})\\
&\leqslant  C\|\phi-\phi'\|_{X^{\tilde\epsilon}}(1+\|\phi\|_{X^{\tilde\epsilon}}^{\rho-1}+\|\phi'\|_{X^{\tilde\epsilon}}^{\rho-1})\\
&\leqslant C\Big\|\begin{bmatrix} \phi\\ \varphi \end{bmatrix}-\begin{bmatrix} \phi'\\ \varphi' \end{bmatrix}\Big\|_{Y_{-1}^{1+\epsilon}}(1+\|\phi\|_{X^{\tilde\epsilon}}^{\rho-1}+\|\phi'\|_{X^{\tilde\epsilon}}^{\rho-1})
\end{split}
\]
for some $C>0$, $\rho>1$ and $\gamma(\epsilon)>0$ with $\rho\tilde\epsilon\leqslant \frac{\gamma(\epsilon)}{3}<\frac13$ $(\tilde\epsilon=\frac{\epsilon}{3})$. 

Hence, by Lemma \ref{asa9345y} we obtain \eqref{asa034mt56} and the proof of the result is complete. \qed

\begin{definition}
Let $\epsilon\geqslant0$, $\tau>0$, $\begin{bmatrix} v_0\\ w_0 \end{bmatrix}\in Y_{-1}^1$. We say that $\begin{bmatrix} v\\ w \end{bmatrix}:[t_0,\tau]\to Y_{-1}^1$ is an $\epsilon-$regular mild solution ($\epsilon-$solution for shor) to \eqref{Eq3a}-\eqref{Eq3b} if $\begin{bmatrix} v\\ w \end{bmatrix}\in C([t_0,\tau],Y_{-1}^1)\cap C((t_0,\tau],Y_{-1}^{1+\epsilon})$, and $\begin{bmatrix} v\\ w \end{bmatrix}(t)$ satisfies
\[
\begin{bmatrix} v(t)\\ w(t) \end{bmatrix}=e^{-\varLambda(t-t_0)}\begin{bmatrix} v_0\\ w_0 \end{bmatrix}+\int_{t_0}^t e^{-\varLambda(t-s)}G(\begin{bmatrix} v(s)\\ w(s)\end{bmatrix})ds.
\]
\end{definition}

Thanks to \cite[Theorem 1]{ACar}, see also \cite[Theorem 3]{CCCC}, we have the following on existence of $\epsilon$-regular solution to \eqref{Eq3a}-\eqref{Eq3b} on certain interval $[0,\tau]$.

\begin{theorem}\label{as024mb84}
Let $\epsilon\geqslant0$, $\tau>0$, $\begin{bmatrix} v_0\\ w_0 \end{bmatrix}\in Y_{-1}^1$.  Then exists a unique $\epsilon$-regular solution to \eqref{Eq3a}-\eqref{Eq3b} on certain interval $[0,\tau]$. This solution satisfies
\[
\begin{bmatrix} v\\ w \end{bmatrix}\in C((0,\tau],Y_{-1}^1)\cap C((0,\tau],Y_{-1}^{1+\theta}),\quad 0\leqslant \theta<\gamma(\epsilon),
\]
and
\[
t^\theta\|\begin{bmatrix} v(t)\\ w(t) \end{bmatrix}\|_{Y_{-1}^{1+\theta}}\to0,\quad\mbox{as}\quad t\searrow0,\quad 0< \theta<\gamma(\epsilon).
\]
\end{theorem}

From this, we establish local well posedness for the Cauchy problem  \eqref{Eq4aaa}-\eqref{Eq4bbb}. We would like to study this problem in the phase space $Z$. To pose the problem in the mentioned space we will need to consider the nonlinear term $F$ as a map with values in the extrapolated space $Z_{-1}$ associated to $\mathbb{A}_\eta$ in $Z$.

\begin{definition}
Let $\epsilon=0$, $\tau>0$, $\Big[\begin{smallmatrix} u_0\\ v_0\\ w_0 \end{smallmatrix}\Big]\in Z^1_{-1}$ and $\Big[\begin{smallmatrix} u\\ v\\ w \end{smallmatrix}\Big](\cdot):[0,\tau]\to Z_{-1}^1$. We say that $\Big[\begin{smallmatrix} u\\ v\\ w \end{smallmatrix}\Big](\cdot)$ is an $\epsilon-$regular solution to \eqref{Eq4aaa}-\eqref{Eq4bbb} on $[0,\tau]$ if and only if 
\begin{itemize}
\item[i)] $\Big[\begin{smallmatrix} u\\ v\\ w \end{smallmatrix}\Big](\cdot)\in C([0,\tau], Z_{-1}^1)\cap C((0,\tau], Z_{-1}^{1+\epsilon})$;
\item[ii)] $\Big[\begin{smallmatrix} u\\ v\\ w \end{smallmatrix}\Big](\cdot)$ satisfies the Cauchy integral formula:
\begin{equation}\label{Semasas_1}
\Big[\begin{smallmatrix} u\\ v\\ w \end{smallmatrix}\Big](t)=e^{-\mathbb{A}_\eta t}\Big[\begin{smallmatrix} u_0\\ v_0\\ w_0 \end{smallmatrix}\Big]+\displaystyle\int_0^t e^{-\mathbb{A}_\eta (t-s)}F\Big(\Big[\begin{smallmatrix} u\\ v\\ w \end{smallmatrix}\Big](s)\Big)ds,\ t\in[0,\tau].
\end{equation}
\end{itemize}
\end{definition}

\begin{theorem}\label{MTheeoo}
Let $\epsilon=0$, $\tau>0$,  $\Big[\begin{smallmatrix} u_0\\ v_0\\ w_0 \end{smallmatrix}\Big]\in  Z^1_{-1}$. Then there exists a unique $\epsilon-$regular solution to \eqref{Eq4aaa}-\eqref{Eq4bbb} on certain interval $[0,\tau]$. In addition, we have 
\[
\Big[\begin{smallmatrix} u\\ v\\ w  \end{smallmatrix}\Big]\in C((0,\tau], Z_{-1}^{1})\cap C((0,\tau], Z_{-1}^{1+\theta}),\quad 0\leqslant \theta<\gamma(\epsilon),
\]
and
\[
t^\theta\|\begin{bmatrix} u(t)\\ v(t)\\ w(t) \end{bmatrix}\|_{Z_{-1}^{1+\theta}}\to0,\quad\mbox{as}\quad t\searrow0,\quad 0< \theta<\gamma(\epsilon).
\]
\end{theorem}

\proof
We also know that, for $\epsilon>0$, $t^\epsilon \left\|e^{-\Lambda t}\left[ \begin{matrix} v_0 \\ w_0 \end{matrix} \right]\right\|_{Y^{\epsilon}}\stackrel{t\to 0}{\longrightarrow}0$. Knowing that $Y^{\epsilon}=X^{\frac13+\epsilon}\times X^\epsilon$, we have that $t^\epsilon \|v(t,v_0,w_0)\|_{X^{\frac13+\epsilon}}\stackrel{t\to 0}{\longrightarrow}0$. Then the integral defining $u(t,u_0,v_0,w_0)$ is convergent and the resulting function is continuous at $t=0$. In fact
\[
\begin{split}
\|u-u_0\|_{X^\frac23} &\leqslant \|e^{-A^{\frac13}t}u_0-u_0\|_{X^\frac23}+\displaystyle\int_0^t \| e^{-A^{\frac13}(t-s)}\|_{\mathcal{L}(X,X^{1-\epsilon})} \|v(s)\|_{X^{\frac13+\epsilon}}ds\\
 &\leqslant \|e^{-A^{\frac13}t}u_0-u_0\|_{X^\frac23}+M\displaystyle\int_0^t (t-s)^{1-\epsilon} s^{-\epsilon} s^\epsilon\|v(s)\|_{X^{\frac13+\epsilon}}ds\\
  &\leqslant \|e^{-A^{\frac13}t}u_0-u_0\|_{X^\frac23}+M\displaystyle\int_0^1 (1-s)^{-1+\epsilon} s^{-\epsilon} ds 
  \ {\displaystyle \sup_{s\in [0,t]}s^\epsilon\|v(s)\|_{X^{\frac13+\epsilon}}}
\end{split}
\]
This ensures the continuity.

Finally, the result follows from \eqref{cv1} and Theorem \ref{as024mb84}.
\qed

\section{Remaks on the case $\eta=1$}\label{S1}

In this section we consider the case $\eta=1$. We note that the initial value problem associated with equation \eqref{E00q1Linear3} as the Cauchy problem in $Z$
\begin{equation}\label{Eq4aadfdfdfabb111}
\dfrac{d}{dt}\Big[\begin{smallmatrix} u\\ v\\ w \end{smallmatrix}\Big]+ \mathbb{B}_{(1)}\Big[\begin{smallmatrix} u\\ v\\ w \end{smallmatrix}\Big]=F_1(\Big[\begin{smallmatrix} u\\ v\\ w \end{smallmatrix}\Big]),\ t>0,
\end{equation}
and
\begin{equation}\label{Eq4bbbfdfdbbbb111}
\Big[\begin{smallmatrix} u\\ v\\ w \end{smallmatrix}\Big](0)=\Big[\begin{smallmatrix} u_0\\ v_0\\ w_0 \end{smallmatrix}\Big],
\end{equation}
where $v=u_t$ and $w=v_t$ and the unbounded linear operator $\mathbb{B}_{(1)}:D(\mathbb{B}_{(1)})\subset Z\to Z$ is defined by
\begin{equation}\label{Eq4aaabb111}
D(\mathbb{B}_{(1)})=X^1\times X^{\frac23}\times X^{\frac13},
\end{equation}
and
\begin{equation}\label{Eq4bbbbbbb111}
\mathbb{B}_{(1)}\Big[\begin{smallmatrix} u\\ v\\ w \end{smallmatrix}\Big]:= \Big[\begin{smallmatrix} A^{\frac13} & -I & 0 \\ 0 & 0 & - I\\ 0 &  A^{\frac23} & 0\end{smallmatrix}\Big] \Big[\begin{smallmatrix} u\\ v\\ w \end{smallmatrix}\Big]=\Big[\begin{smallmatrix} A^{\frac13}u-v\\ -w\\  A^{\frac23}v  \end{smallmatrix}\Big],\ \forall \Big[\begin{smallmatrix} u\\ v\\ w \end{smallmatrix}\Big]\in X^1\times X^{\frac23}\times X^{\frac13}.
\end{equation}

The nonlinearity $F_1$ given by \eqref{DefF}, where $f:D(A^{\frac23})\subset X\to X$ is a Lipschitz continuous function on bounded sets.

Thanks to Lemma \ref{MleNewddsd}, in particular, we have 
\begin{lemma}\label{MleNewddsd=1}
Let $\mathbb{B}_{(1)}$ be the unbounded linear operator defined in \eqref{Eq4aaabb111}-\eqref{Eq4bbbbbbb111}. The following conditions hold.
\begin{itemize}
\item[$i)$] The unbounded linear operator $\mathbb{B}_{(1)}$ is closed and densely defined;
\item[$ii)$] Zero belongs to the resolvent set $\rho(\mathbb{B}_{(1)})$; namely, the resolvent operator of $\mathbb{B}_{(1)}$ is the bounded linear operator $\mathbb{B}_{(1)}^{-1}:Z\to Z$ given by
\begin{equation}\label{sdsfew5436}
\mathbb{B}_{(1)}^{-1}=\Big[\begin{smallmatrix}  A^{-\frac13} & 0 & A^{-1}\\ 0& 0& A^{-\frac23}  \\ 0 & -I & 0 \end{smallmatrix}\Big].
\end{equation}
Moreover, $\mathbb{B}_{(1)}$ has compact resolvent;
\item[$iii)$] The spectrum set of $-\mathbb{B}_{1}$, $\sigma(-\mathbb{B}_{(1)})$, is given by
\[
\sigma(-\mathbb{B}_{(1)})=\{\lambda \in \mathbb{C}: \lambda \in \sigma(-A^{\frac13})\}\cup\{\lambda i \in \mathbb{C}: \lambda \in \sigma(-A^{\frac13})\}\cup\{-\lambda i \in \mathbb{C}: \lambda \in \sigma(-A^{\frac13})\},
\]
where $\sigma(-A^{\frac13})$ denote the spectrum set of  $-A^{\frac13}$.
\end{itemize}

	\begin{figure}[H]
		\begin{center}
			\begin{tikzpicture}
			\draw[-stealth'] (-5,0) -- (2,0) node[below] {$\scriptstyle {\rm Re}$};
			\draw[-stealth'] (0,-3.15) -- (0,3.15) node[above] {\color{black}$\scriptstyle{\rm Im}$};
			\draw (0,0) -- (-3, 0);
\foreach \i in {0.1,0.14,...,1.2}{\EXP[\i]{3}{\sol}\fill [color=blue] (-3*\sol,0) circle (1.5pt);}
\foreach \i in {0.1,0.14,...,1.2}{\EXP[\i]{3}{\sol}\fill [color=blue] (0,2*\sol) circle (1.5pt);}
\foreach \i in {0.1,0.14,...,1.2}{\EXP[\i]{3}{\sol}\fill [color=blue] (0,-2*\sol) circle (1.5pt);}
\node at (-2.5,0.4) {$\sigma(-A^{\frac13})$};
\node at (-1,2) {$i\sigma(-A^{\frac13})$};
\node at (-1.15,-2) {$-i\sigma(-A^{\frac13})$};
			\end{tikzpicture}
		\end{center}
\caption{Semi-lines contained the eigenvalues of $-\mathbb{B}_{1}$.}\label{figadadaf1=1}
	\end{figure}
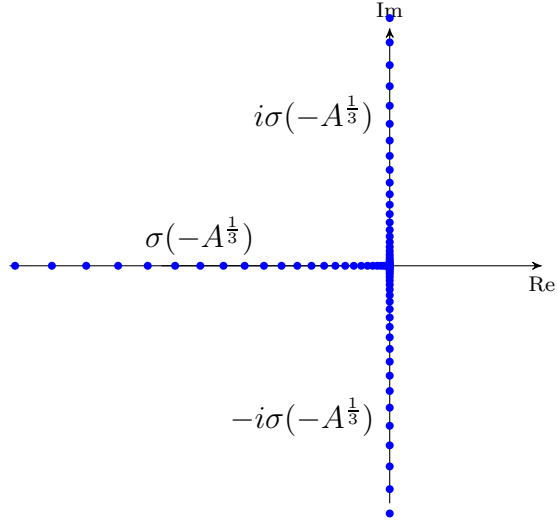
\end{lemma}

\proof
The  proof of $i)$  easily follows from \eqref{Eq4bbbbbbb111}. For the proof of $ii)$ the result easily follows from \eqref{sdsfew5436}. Finally, for the proof of $iii)$ the result easily follows from \eqref{Eq4bbbbbbb111}.
\qed

\begin{theorem}\label{tdissipative}
Let $\mathbb{B}_{(1)}$ be the unbounded linear operator defined in \eqref{Eq4aaabb111}-\eqref{Eq4bbbbbbb111}. The unbounded linear operator $-\mathbb{B}_{(1)}$ is not a dissipative operator on the state space $Z$.
\end{theorem}

\proof
Indeed, let $u\in X^1$ be such that  $u\neq0$ and let $\Big[\begin{smallmatrix} u\\ 2A^{\frac13}u\\ 0 \end{smallmatrix}\Big]\in Z^1$. Note that
\[
\Big\langle \mathbb{B}_{(1)}\Big[\begin{smallmatrix} u\\ 2A^{\frac13}u\\ 0\end{smallmatrix}\Big],\Big[\begin{smallmatrix} u\\ 2A^{\frac13}u\\ 0 \end{smallmatrix}\Big]\Big\rangle=-\langle A^{\frac13}u,u\rangle_{X^{\frac23}}
\]
and consequently
\[
Re\Big\langle -\mathbb{B}_{(1)}\Big[\begin{smallmatrix} u\\ 2A^{\frac13}u\\ 0\end{smallmatrix}\Big],\Big[\begin{smallmatrix} u\\ 2A^{\frac13}u\\ 0 \end{smallmatrix}\Big]\Big\rangle>0
\]
and the prove is complete. 
\qed

Theorem\ref{tdissipative} ensures that the linear operator $-\mathbb{B}_{(1)}$ is not an infinitesimal generator of a specific type of strongly continuous semigroup in $Z$; that is, the strongly continuous semigroup of contractions in $Z$.

\section{Moore-Gibson-Thompson-type equations}\label{FinalSec}

In this subsection we present  boundary-initial value problems associated with a Moore-Gibson-Thompson equation with fractional damped and  strongly damped linear wave equation, where  our results from previous section can be applied. Namely, let $\Omega\subset\mathbb{R}^N$, $N\geqslant3$, be a bounded smooth domain, and the initial-boundary value problems
\begin{equation}\label{Prob_1App}
\begin{cases}
u_{ttt}  -\Delta u + \eta  (-\Delta)^{\frac13} u_{tt} + \eta  (-\Delta)^{\frac23} u_t=f(u),&\ t>0,\ x\in\Omega,\\
u(0,x)=\varphi(x),\ u_t(0,x)=\xi(x),\ u_{tt}(0,x)=\psi(x),&\ x\in\Omega,\\
u(t,x)=0,&\ t\geqslant 0,\ x\in\partial\Omega,
\end{cases}
\end{equation}
where $\eta>1$. 

The nonlinearity $f:\mathbb{R}\to\mathbb{R}$ in \eqref{Prob_1App} is a continuously differentiable function satisfying for some $1<\rho< \frac{3N+4}{3N-8}$ the growth condition
\begin{equation}\label{CondfDeriv}
|f'(s)|\leqslant C(1+|s|^{\rho-1}).
\end{equation}

Here we consider $X=L^2(\Omega)$ and the negative Laplacian operator
\[
A u=-\Delta u,
\]
with domain
\[
D(A)=H^2(\Omega)\cap H^1_0(\Omega)
\]
which is a sectorial operator and it bounded imaginary powers, and consequently the spaces $X^\alpha$, $\alpha\in[0,1]$, are characterized with the aid of complex interpolation as
\[
X^\alpha=[H^2(\Omega)\cap H^1_0(\Omega),L^2(\Omega)]_\alpha=H^{2\alpha}_{\{I\}}(\Omega)
\]
and
\[
X^{-\alpha}=(H^{2\alpha}_{\{I\}}(\Omega))'
\]
where $[\cdot,\cdot]_\alpha$ denotes the complex interpolation function (see \cite{A} and \cite{JDu}). In particular $X=X^0=L^2(\Omega)$, $X^{\frac12}=H_0^1(\Omega)$,  $X^{-\frac12}=(H_0^1(\Omega))'$  and $X^1=H^2(\Omega)\cap H^1_0(\Omega)$. 

With this set-up we will consider problem \eqref{Prob_1App} in the form \eqref{Eq4aadfdfdfabb}-\eqref{Eq4bbbfdfdbbbb} with $u_0=\varphi, v_0=A^\frac13\varphi+\xi$, and $w_0=A^\frac13\xi+\psi$.

Let $F: Z_{-1}^1\to Z_{-1}^\alpha$, $\alpha\geqslant0$, be a locally Lipschitz continuous map, as well as in  \eqref{DefF}. Recall that a mild solution of  \eqref{Eq4aadfdfdfabb}-\eqref{Eq4bbbfdfdbbbb} on $[0,\tau]$ is a function $\Big[\begin{smallmatrix} u\\ v\\ w \end{smallmatrix}\Big]\Big(\cdot,\Big[\begin{smallmatrix} u_0\\ v_0\\ w_0 \end{smallmatrix}\Big]\Big)\in C([0,\tau],Z_{-1}^1)$ which satisfies 
\[
\Big[\begin{smallmatrix} u\\ v\\ w \end{smallmatrix}\Big]\Big(t,\Big[\begin{smallmatrix} u_0\\ v_0\\ w_0 \end{smallmatrix}\Big]\Big)=e^{-\mathbb{B}_\eta t}\Big[\begin{smallmatrix} u_0\\ v_0\\ w_0 \end{smallmatrix}\Big]+\displaystyle\int_0^t e^{-\mathbb{B}_\eta (t-s)}F\Big(\Big[\begin{smallmatrix} u\\ v\\ w \end{smallmatrix}\Big]\Big(s,\Big[\begin{smallmatrix} u_0\\ v_0\\ w_0 \end{smallmatrix}\Big]\Big)\Big)ds,
\]
for $t\in[0,\tau]$. We say that \eqref{Eq4aadfdfdfabb}-\eqref{Eq4bbbfdfdbbbb} is locally well posed in $Z_{-1}^1$ is for  any $\Big[\begin{smallmatrix} u_0\\ v_0\\ w_0 \end{smallmatrix}\Big]\in Z_{-1}^1$ there is a unique mild solution 
\[
t\mapsto \Big[\begin{smallmatrix} u\\ v\\ w \end{smallmatrix}\Big]\Big(t,\Big[\begin{smallmatrix} u_0\\ v_0\\ w_0 \end{smallmatrix}\Big]\Big)
\]
of \eqref{Eq4aadfdfdfabb}-\eqref{Eq4bbbfdfdbbbb} defined on a maximal interval of existence $[0,t_{u_,v_0,w_0})$ and depending continuously on the initial data $\Big[\begin{smallmatrix} u_0\\ v_0\\ w_0 \end{smallmatrix}\Big]$.

As a consequence of the Sobolev embeddings we obtain the following result cf. \cite[Proposition 1.3.8]{ChD}.

\begin{proposition}
Let $\Omega\subset\mathbb{R}^N$ be a bounded domain of class $C^m$ and $(A,D(A))$ be a sectorial operator in $L^p(\Omega)$, $1<p<\infty$, with $D(A)\subset W^{2m,p}(\Omega)$ for some $m\geqslant1$. Then for $\alpha\in[0,1]$ the following inclusion holds.
\[
X^\alpha\subset W^{s,q}(\Omega)
\]
if $2m\alpha-\frac{N}{p}\geqslant s-\frac{N}{q}$, $1<p\leqslant q<\infty$, $s\geqslant0$.
\end{proposition}

From we have

\begin{theorem}
The problem \eqref{Eq4aadfdfdfabb}-\eqref{Eq4bbbfdfdbbbb} with $\eta>1$ is locally posed in $Z_{-1}^1$ whenever $f$ satisfies \eqref{CondfDeriv} for some $1<\rho< \frac{3N+4}{3N-8}$
\end{theorem}

\proof
The map $F$ defined as in \eqref{DefF} is Lipschitz continuous on bounded sets from $Z_{-1}^1$ into $Z_{-1}^{\sigma}=X^{\frac{1+\sigma}{3}}\times X^{\frac{\sigma}{3}}\times X^{\frac{-1+\sigma}{3}}$ whenever $0<\sigma\leqslant \tilde{\sigma}$, and $\tilde{\sigma}=\min\{1, (\rho-1)(2-\frac{3}{4}N)+1 \}$ ($X^{\frac{2}{3}}\hookrightarrow L^{\frac{6N(\rho-1)}{8+4(2\sigma)}}(\Omega)$). Indeed, if $B$ is a bounded  subset of $Z_{-1}^1$ and $\Big[\begin{smallmatrix} u_1\\ v_1\\ w_1 \end{smallmatrix}\Big],\Big[\begin{smallmatrix} u_2\\ v_2\\ w_2 \end{smallmatrix}\Big] \in B$, we have
\[
\Big\| F\Big( \Big[\begin{smallmatrix} u_1\\ v_1\\ w_1 \end{smallmatrix}\Big]\Big) - F\Big(\Big[\begin{smallmatrix} u_2\\ v_2\\ w_2 \end{smallmatrix}\Big] \Big)  \Big\|_{Z_{-1}^\sigma}\leqslant c_1\|f(u_1)-f(u_2)\|_{X^{\frac{-1+\sigma}{3}}}.
\]
Since $L^{\frac{6N}{3N+4(1-\sigma)}}(\Omega) \hookrightarrow H^{\frac{2(1-\sigma)}{3}}_{\{I\}}(\Omega)=X^{\frac{-1+\sigma}{3}}$ we obtain  
\[
\Big\| F\Big( \Big[\begin{smallmatrix} u_1\\ v_1\\ w_1 \end{smallmatrix}\Big]\Big) - F\Big(\Big[\begin{smallmatrix} u_2\\ v_2\\ w_2 \end{smallmatrix}\Big] \Big)  \Big\|_{Z_{-1}^\sigma}\leqslant c_2\|f(u_1)-f(u_2)\|_{L^{\frac{6N}{3N+4(1-\sigma)}}(\Omega)}
\]
and thanks to \eqref{CondfDeriv} there exists $C>0$ such that 
\[
\forall s_1,s_2\in\mathbb{R},\quad | f(s_1)-f(s_2)|\leqslant C|s_1-s_2|(1+|s_1|^{\rho-1}+|s_2|^{\rho-1})
\]
and consequently 
\[
\begin{split}
\Big\| F\Big( \Big[\begin{smallmatrix} u_1\\ v_1\\ w_1 \end{smallmatrix}\Big]\Big) - F\Big(\Big[\begin{smallmatrix} u_2\\ v_2\\ w_2 \end{smallmatrix}\Big] \Big)  \Big\|_{Z_{-1}^\sigma}&\leqslant c_3\|u_1-u_2\|_{X^{\frac23}}\Big(1+\|u_1\|^{\rho-1}_{L^{\frac{6N(\rho-1)}{8+4(1-\sigma)}}(\Omega)}+\|u_2\|^{\rho-1}_{L^{\frac{6N(\rho-1)}{8+4(1-\sigma)}}(\Omega)}\Big)\\
&\leqslant c_4\Big\|  \Big[\begin{smallmatrix} u_1\\ v_1\\ w_1 \end{smallmatrix}\Big] -\Big[\begin{smallmatrix} u_2\\ v_2\\ w_2 \end{smallmatrix}\Big]   \Big\|_{Z_{-1}^1}.
\end{split}
\]

The proof now follows from \cite{He}.
\qed

\end{document}